\documentclass[11pt]{amsart}

\usepackage{amscd,amsfonts,amssymb}
\usepackage{tikz}
\usepackage{pgfplots}
\usetikzlibrary{automata,positioning,calc,trees}
\usetikzlibrary{intersections,pgfplots.fillbetween}
\usepackage{graphicx,tikz}
\usepackage{pgf,tikz,pgfplots}
\usepackage{mathrsfs}
\usetikzlibrary{arrows}
\usepackage{enumerate}
\usepackage[shortlabels]{enumitem}
\usepackage{mathrsfs}
\usepackage{amssymb,amsmath,amsthm,color}
\usepackage{caption,subcaption}
\usepackage[alphabetic,bibtex-style]{amsrefs}
\usepackage{hyperref}
\usepackage{url}
\usepackage{setspace}
\usepackage{float}

\BibSpec{article}{%
	+{}  {\PrintAuthors}                {author}
	+{,} { \textit}                     {title}
	+{.} { }                            {part}
	+{:} { \textit}                     {subtitle}
	+{,} { \PrintContributions}         {contribution}
	+{.} { \PrintPartials}              {partial}
	+{,} { }                            {journal}
	+{}  { \textbf}                     {volume}
	+{}  { \PrintDatePV}                {date}
	+{,} { \issuetext}                  {number}
	+{,} { \eprintpages}                {pages}
	+{,} { }                            {status}
	+{,} { \url}                        {url}    
	+{,} { \PrintDOI}                   {doi}
	+{,} { available at \eprint}        {eprint}
	+{}  { \parenthesize}               {language}
	+{}  { \PrintTranslation}           {translation}
	+{;} { \PrintReprint}               {reprint}
	+{.} { }                            {note}
	+{.} {}                             {transition}
}

\date{\today}
\allowdisplaybreaks
\newcommand{\R}{{\mathbb R}}       
\newcommand{\N}{{\mathbb N}}

\newcommand{\Z}{{\mathbb Z}}       

\newtheorem{theorem}{Theorem}[section]
\newtheorem{lemma}[theorem]{Lemma}
\newtheorem{remark}[theorem]{Remark}

\newtheorem{proposition}[theorem]{Proposition}

\newtheorem*{theorem*}{Theorem}

\theoremstyle{definition}

\numberwithin{equation}{section}

\begin{document}
	
\title{On the bilinear cone multiplier}

\author[Shrivastava]{Saurabh Shrivastava}
\address{Saurabh Shrivastava\\
	Department of Mathematics\\
	IISER Bhopal\\
	Bhopal-462066, India}
\email{saurabhk@iiserb.ac.in}
\author[Shuin]{Kalachand Shuin}
\address{Kalachand Shuin\\
	Department of Mathematics\\
	Indian Institute of Science\\
	Bengaluru-560012, India.}
\email{kalachands@iisc.ac.in, shuin.k.c@gmail.com}
	\begin{abstract}
	For $f,g \in \mathcal S(\R^n), n\geq 3$, consider the bilinear cone multiplier operator defined by 
	\[\mathcal{T}^{\lambda}_{R}(f,g)(x):=\int_{\mathbb{R}^{2n}}m^{\lambda}\left(\frac{\xi'}{R\xi_n},\frac{\eta'}{R\eta_n}\right)\hat{f}(\xi)\hat{g}(\eta)e^{2\pi\iota x\cdot(\xi+\eta)}~d\xi d\eta,\]
	where $\lambda>0, R>0$ and \[m^{\lambda}\left(\frac{\xi'}{R\xi_n},\frac{\eta'}{R\eta_n}\right)=\Big(1-\frac{|\xi'|^2}{R^2\xi^2_n}-\frac{|\eta'|^2}{R^2\eta^2_n}\Big)^{\lambda}_{+}\varphi(\xi_n)\varphi(\eta_n),\]
	$(\xi',\xi_n), (\eta',\eta_n)\in\mathbb{R}^{n-1}\times \mathbb{R}$ and  $\varphi\in C_{c}^{\infty}([\frac{1}{2},2])$.
	
	We investigate the problem of pointwise almost everywhere convergence of $\mathcal{T}^{\lambda}_{R}(f,g)(x)$ as $R\rightarrow \infty$ for $(f,g)\in L^{p_1}\times L^{p_2}$ for a wide range of exponents $p_1, p_2$ satisfying the H\"{o}lder relation $\frac{1}{p_1}+\frac{1}{p_2}=\frac{1}{p}$. This assertion is proved by establishing suitable weighted $L^{2}\times L^{2}\rightarrow L^{1}$--estimates of the maximal bilinear cone multiplier operator 
	\[\mathcal{T}^{\lambda}_{*}(f,g)(x):=\sup_{R>0}|\mathcal{T}^{\lambda}_{R}(f,g)(x)|.\]
	\end{abstract}

	\maketitle
	
\section{introduction}
\subsection{Bochner--Riesz means}
The convergence of the $N$th spherical partial sum operator \[S_{N}f(x):=\sum_{|\textbf{n}|\leq N}  \hat{f}(\textbf{n}) e^{2\pi\iota x\cdot\textbf{n}},\]
where  $f\in C(\mathbb T^n), \textbf{n}\in\mathbb{Z}^n, N\in\mathbb{N}$, motivates the study of maximal Fourier multiplier operators. In general, one would like to consider the Fourier multiplier operators of the form 
\[S_{\phi,N}f(x):=\sum_{\textbf{n}\in \Z^n} \Phi(\textbf{n}, N) \hat{f}(\textbf{n}) e^{2\pi\iota x\cdot\textbf{n}},\]
where $\Phi(\cdot, N)$ is a suitable  bounded function for each $N$. The uniform $L^p$--estimates of $S_{\Phi,N}f$ and the associated maximal function $\sup\limits_{N\in \N}|S_{\Phi,N}f|$ imply the convergence of $S_{\Phi,N}$ in $L^p$--norm and pointwise almost everywhere respectively. In the recent past, there have been many interesting developments in this area providing deeper insights into the relationship between analytic and geometric properties of multipliers $\Phi(\textbf{n},N)$.

The well-known transference principle for Fourier multipliers allows us to address these questions via the study of analogous problems for the corresponding Fourier multipliers defined on $\R^n.$  More precisely, the continuous analogue of the partial sum Fourier multiplier operator is defined by 
\[S_{\Phi,R}f(x):=\int_{\mathbb{R}^{n}}\Phi(\xi, R)\hat{f}(\xi) e^{2\pi\iota x\cdot\xi}~d\xi, x\in \R^n, R>0,\]
where $\Phi(\cdot, R)$ is a suitable  bounded measurable function on $\R^n$ for each $R>0$. The case of spherical partial sum operator of Fourier series corresponds to the study of the Fourier multiplier operator $S_{\Phi,R}f$, where $\Phi(\xi,R)=\chi_{B(0,R)}(\xi)$, and $B(0,R)=\{x\in \R: |x|<R\}$ is the euclidean ball with its center at the origin and radius $R$. In this case we denote $S_{\Phi,R}=S_{R}$. It is a much celebrated result due to C. Fefferman~\cite{Fefferman} that the $L^p$--estimates 
\[\|S_{R}f\|_p\lesssim \|f\|_p,\]
do not hold when $p\neq 2$ and $n\geq 2.$ 
In view of this result, various summability methods were introduced to study the convergence properties of Fourier series and Fourier multipliers. In particular, the Bochner--Riesz means are of extremely importance in the subject.   
Bochner and Riesz introduced the new partial sum multiplier operator \[S^{\lambda}_{N}f(x):=\sum_{|\textbf{n}|\leq N} \Big(1-\frac{|\textbf{n}|}{N}\Big)^{\lambda} \hat{f}(\textbf{n}) e^{2\pi\iota x\cdot\textbf{n}},\]
for $\lambda\geq0$. As mentioned above, using the transference principle \cite{deLeeuw} it suffices to study $L^p$-boundedness of the continuous analogue of Bochner--Riesz  multiplier operator, defined by 

\[B^{\lambda}_{R}f(x):=\int_{\mathbb{R}^{n}}\Big(1-\frac{|\xi|^{2}}{R^{2}}\Big)^{\lambda}_{+}\hat{f}(\xi) e^{2\pi\iota x\cdot\xi}~d\xi, \lambda\geq0~~\text{and}~ f\in\mathcal{S}(\mathbb{R}^{n}), \]
where the notation $r_{+}=\max\{r,0\}$ and $R>0$.

The Bochner--Riesz conjecture states that $B^{\lambda}_R$ is bounded from $L^{p}(\mathbb{R}^n)$ to itself with uniform bounds with respect to the parameter $R$ if, and only if $\lambda>\max\{n|\frac{1}{p}-\frac{1}{2}|-\frac{1}{2},0\}$.
For $n=2$, the conjecture was resolved positively by Carleson and  Sj\"{o}lin \cite{CS} and Carbery \cite{Carbery}.
For the higher dimensional case   (i.e. $n\geq3$), there have been many interesting developments in this direction in the recent past, however the conjecture remains unresolved till date.  The reader is referred to \cite{ LRS, LW, Lee0, Lee1, Christ, GW, Tao} and references therein for more details.

The study of pointwise almost everywhere convergence of the Bochner--Riesz means leads to investigation of  suitable $L^p$--estimates of the maximal Bochner--Riesz mean defined by  \[B^{\lambda}_{*}f(x):=\sup_{R>0}|B^{\lambda}_{R}f(x)|.\] 
Carbery, Rubio de Francia and Vega proved the following result.
\begin{theorem*}[\cite{CarberyRubioVega}]
	Let $\lambda>0$ and $n\geq2$. Then for all $f\in L^{p}(\mathbb{R}^n)$ with $2\leq p<\frac{2n}{n-1-2\lambda}$, we have 
	\[\lim_{R\to\infty} B^{\lambda}_{R}f(x)=f(x) ~\text{a.e.}~x.\]
\end{theorem*}
They deduced the result stated above by proving the following weighted $L^2$-estimate of $B^{\lambda}_{*}$.
\begin{theorem*}[\cite{CarberyRubioVega}]
	Let $\lambda>0$ and $0\leq \alpha<1+2\lambda\leq n$. Then 
	\[\int_{\mathbb{R}^n} |B^{\lambda}_{*}f(x)|^2 |x|^{-\alpha}~dx\leq C_{\alpha,\lambda}\int_{\mathbb{R}^n} |f(x)|^2 |x|^{-\alpha}~dx.\]
\end{theorem*}
We refer to some of the recent papers \cites{Jongchon, Xiaochun} for $L^{p}$--estimates of the maximal Bochner--Riesz mean for $1<p<2$.

\subsection{Cone multipliers}
The problem of cone multiplier operator was introduced by Stein~\cite{SteinSomeproblemsinharmonicanalysis}. Let $n\geq 3$ and write $\xi=(\xi',\xi_n)\in \R^n$, where $\xi'\in \R^{n-1}$ and $\xi_n\in \R$. For  $\lambda>0$ and  $R>0$, the cone multiplier operator is defined by  
\[T_{R}^{\lambda}f(x):=\int_{\mathbb{R}^n}\Big(1-\frac{|\xi'|^2}{R^2\xi^2_{n}}\Big)^{\lambda}_{+}\varphi(\xi_n)\hat{f}(\xi) e^{2\pi\iota x\cdot\xi}~d\xi,\]
where $\varphi\in C_{c}^{\infty}([1/2,2])$ is a non-negative function. 

Observe that for fixed $\xi_n$, the cone multiplier $(1-\frac{|\xi'|^2}{R^2\xi^2_{n}})^{\lambda}_{+}$ in $\R^n$ reduces to Bochner--Riesz multiplier in dimension $n-1$. Indeed, it is conjectured that the range of $p$ for $L^p$-boundedness of the cone multiplier operator $T_{R}^{\lambda}$ in $\R^n$ is precisely the same as that of the Bochner--Riesz multiplier operator $B^{\lambda}_{R}$ in $\R^{n-1}$. 
Mockenhaupt\cite{Mockenhaupt} reduced the $L^p$-boundedness of the operator $T_{R}^{\lambda}$ to the study of reverse square function. He proved an $L^4$-estimate of the operator when $\lambda>\frac{1}{8}.$ Later, Bourgain~\cite{Bourgain} introduced the bilinear method and obtained improved $L^p$--estimates of the operator  $T_{R}^{\lambda}.$ These methods have been generalized further to address the local smoothing conjecture and  Fourier restriction problem. There have been many interesting developments in the direction of cone multiplier operator, for example, see the papers \cites{Garrigos,Heo,HHY,LeeVargas, CHLY,GuthWangZhangAsharpSquareFunctionEstimate} for some important developments on the cone multiplier operator and related topics. In particular, the breakthrough paper by Guth, Wang and Zhang~\cite{GuthWangZhangAsharpSquareFunctionEstimate} established sharp estimates of the reverse square function. Consequently, they resolved the local smoothing conjecture and cone multiplier conjecture in dimension $n=3$.

Recently, Chen, He, Li, and Yan~\cite{CHLY} investigated the pointwise almost everywhere convergence of the cone multiplier. In order to do so they employed the approach of \cite{CarberyRubioVega} and proved weighted $L^2$-estimates of the maximal cone multiplier operator, defined by   
\[T^{\lambda}_{*}f(x)=\sup_{R>0}|T^{\lambda}_{R}f(x)|.\]

They proved the following result.
\begin{theorem}\label{Chenelal} \cite{CHLY}
Let $n \geq 3, \alpha \in[0, n-1), \beta \in[0,1)$, and $\lambda>\max \left\{\frac{\alpha+\beta-1}{2}, 0\right\}$.
Then, the following weighted $L^2$-estimate holds 
$$
\left\|T_*^\lambda(f)\right\|_{L^2\left(w_{\alpha, \beta}\right)} \lesssim \|f\|_{L^2\left(w_{\alpha, \beta}\right)}, f\in \mathcal S(\R^n), 
$$
where $w_{\alpha, \beta}(x)=|x'|^{-\alpha}|x_n|^{-\beta}, x=(x',x_n)\in \R^{n-1}\times \R.$
\end{theorem} 
Consequently, they proved the following pointwise almost everywhere result for the cone multiplier. 
\begin{theorem*}[Theorem 1.2, \cite{CHLY}]
	Let $n\geq3$, $p\in [2,\infty)$, and $\lambda>\max\{n|1/p-1/2|-1/2,0\}$. Then for $f\in L^{p}(\mathbb{R}^n)$ with $supp(\hat{f})\subset \{\xi\in\mathbb{R}^n: a<|\xi_n|<b\}, 0<a<b<\infty$, the following convergence holds 
	\begin{eqnarray*}\label{maximalcone}
		\lim_{R\to\infty} T^{\lambda}_{R}f(x)=f(x),~~\text{a.e.}~x \in \R^n.
	\end{eqnarray*}
\end{theorem*}
\subsection{Bilinear Bochner--Riesz means}
In the past two decades the study of bilinear Fourier multipliers has emerged as an active research area in harmonic analysis. In particular, there have been some interesting developments in the direction of the bilinear Bochner--Riesz means. These operators are defined in the following manner. 

Let $\lambda\geq0$ and $R>0$.   The bilinear Bochner--Riesz mean of order $\lambda$ is defined by  
\[\mathcal {B}_{R}^{\lambda}(f,g)(x):=\int_{\mathbb{R}^{2n}}\Big(1-\frac{\vert \xi\vert^{2}+\vert \eta\vert^{2}}{R^{2}}\Big)^{\lambda}_{+}\hat{f}(\xi)\hat{g}(\eta)e^{2\pi \iota x\cdot(\xi+\eta)}d\xi d\eta, f,g\in \mathcal S(\R^n).\]

The study of $L^{p_1}\times L^{p_2}\to L^{p}$--estimates of the operator $\mathcal {B}_{R}^{0},$ known as the bilinear disc multiplier operator, in dimension $n=1$ was initiated by Grafakos and Li~\cite{GrafakosLiTheDiscBilinearMultiplier}. Later, Bernicot, Grafakos, Song and Yang~\cite{BGSY} investigated the $L^{p_1}\times L^{p_2}\to L^{p}$-estimates of the operator $\mathcal {B}_{R}^{\lambda}$ for $\lambda>0$ and $n\geq 1.$ They obtained partial results for the operator under consideration. Since then, many interesting results concerning the bilinear Bochner--Riesz means and related objects have been proved in the recent past. We discuss some of them here. Jeong, Lee and Vargas~\cite{JLV} extended the $L^{p_1}\times L^{p_2}\to L^{p}$--estimates of the operator $\mathcal {B}_{R}^{\lambda}$ for an improved range of exponents $p_1, p_2, p$ and $\lambda>0$. They reduced the study of bilinear Bochner--Riesz means to a product of discretized square functions associated with linear Bochner--Riesz means. Later, Jeong and Lee~\cite{JL} generalized the method of \cite{JLV} further to obtain partial results for the maximal bilinear Bochner--Riesz means, defined by 
\[\mathcal{B}^{\lambda}_{*}(f,g)=\sup_{R>0}\left|\mathcal {B}_{R}^{\lambda}(f,g)(x)\right|.\]
Very recently, Jotsaroop and Shrivastava \cite{JS} introduced a new method to address the problem of $L^p$-estimate of bilinear Bochner--Riesz means and associated maximal function. Their method relies on decomposing the bilinear Bochner--Riesz multiplier into an average of product of shifted Bochner--Riesz means. This allows them to obtain improved $L^p$--estimates of the maximal bilinear Bochner--Riesz operator. This method has been further generalized to study bilinear analogue of Stein's square function for Bochner--Riesz means and multilinear Bochner--Riesz means. We refer the reader to ~\cites{CS,JotsaroopShrivastavaShuinWeightedBilinearBochnerRiesz,ChoudharyJotsaroopShrivastavaSquareFunction,Shuin,HeLiZhengPointwiseConvergenceMultilinearBochnerRiesz, BhojakChoudharyShrivastavaConvexDomains, JotsaroopShrivastavaShuinWeightedBilinearBochnerRiesz} and references therein for some of these developments and related works. 

\section{Main results for bilinear cone multiplier}
The aforementioned papers on cone multiplier problem and recent developments in the context of bilinear Bochner--Riesz means motivate us to consider the bilinear analogue of cone multiplier. In this paper, we introduce these objects and initiate the study of pointwise almost everywhere convergence. 

Let $\lambda>0$ and $R>0$ be two parameters. The bilinear analogue of cone multiplier operator in $\R^n, n\geq 3$ is defined by 
\begin{align*}
	\mathcal{T}^{\lambda}_{R}(f,g)(x):=\int_{\mathbb{R}^{2n}}m^{\lambda}(\frac{\xi'}{R\xi_n},\frac{\eta'}{R\eta_n})\hat{f}(\xi)\hat{g}(\eta)e^{2\pi\iota x\cdot(\xi+\eta)}~d\xi d\eta,
\end{align*}
where $m^{\lambda}(\frac{\xi'}{R\xi_n},\frac{\eta'}{R\eta_n})=\Big(1-\frac{|\xi'|^2}{R^2\xi^2_n}-\frac{|\eta'|^2}{R^2\eta^2_n}\Big)^{\lambda}_{+}\varphi(\xi_n)\varphi(\eta_n)$ for $\varphi \in C_{c}^{\infty}([1/2,2]).$ 
\begin{remark}
	The transference principle for bilinear multipliers could be applied to deduce necessary conditions for $L^{p_1}(\R^n)\times L^{p_2}(\R^n)\rightarrow L^{p}(\R^n)$--boundedness of the bilinear cone multiplier via transferring it to the boundedness of the bilinear Bochner--Riesz mean of same index in lower dimension $\R^{n-1},$ see \cite{Shobha, Oscar, Grafakosclassical} for relevant results on tranference principle in the context of bilinear multipliers. However, the necessary conditions for the bilinear Bochner--Riesz means are not well understood. Therefore, we do not pursue this direction here.
\end{remark}
We prove the following result concerning the pointwise almost everywhere convergence of $\mathcal{T}^{\lambda}_{R}(f,g)(x)$.
\begin{theorem}\label{mainresult}
	Let $n\geq3$ and $2\leq p_1,p_2<\infty$ with $1/p=1/p_1+1/p_2$. Then for $f\in L^{p_1}(\mathbb{R}^{n})$, $g\in L^{p_2}(\mathbb{R}^n)$  with 
	\[supp(\hat{f}),supp(\hat{g})\subset \{\xi\in\mathbb{R}^n: \frac{1}{2}<|\xi_n|<2\}.\]
The following pointwise almost convergence holds  
	\[\lim_{R\to\infty}\mathcal{T}^{\lambda}_{R}(f,g)(x)=f(x)g(x)~~~\text{for a.e. }x\in\mathbb{R}^{n},\]
	whenever $\lambda>\max\{n(1-\frac{1}{p})-1,n(\frac{1}{2}-\frac{1}{p_1})-\frac{1}{2},n(\frac{1}{2}-\frac{1}{p_2})-\frac{1}{2},0\}$.
\end{theorem}
\begin{remark}
	Observe that $\lambda>0$, the result above provides almost everywhere convergence of bilinear cone multipliers $\mathcal{T}^{\lambda}_{R}(f,g)(x)$ for $(f,g)\in  L^{p_1}(\R^n)\times L^{p_2}(\R^n)$ for $p_1, p_2\in [2,\frac{2n}{n-1}]$. 
\end{remark}

The proof of Theorem~\ref{mainresult} is developed following the approach provided in~\cite{CarberyRubioVega} and \cite{CHLY}. Their approach reduces the job to proving suitable weighted $L^2\times L^2\to L^1$-estimate, referred to as weighted $(2,2,1)$ estimate, of the maximal bilinear cone multiplier operator defined by 
\begin{align*}
	 \mathcal{T}^{\lambda}_{*}(f,g)(x):=\sup_{R>0}| \mathcal{T}^{\lambda}_{R}(f,g)(x)|.
\end{align*}
We require some notation in order to state this reduction. 

Let $x=(x',x_n)\in \mathbb{R}^{n-1}\times \mathbb{R}$ and $\alpha,\beta\in \mathbb{R}$. Define $\omega_{\alpha,\beta}(x)=|x'|^{-\alpha}|x_n|^{-\beta}$. It is easy to verify that $\omega_{\alpha,\beta}$ is an Muckenhoupt weight in the class $ A_{p}(\mathbb{R}^{n-1}\times \mathbb{R}),~1<p<\infty$ if, and only if $1-n<-\alpha<(n-1)(p-1), -1<-\beta<p-1$. We refer to~\cite{DuoanBook} for an introduction to Muckenhoupt $A_p$ weights. 

Let $\alpha\in [0,n-1), \beta\in [0,1)$ and $0<\delta\ll 1$. Define the constant \[A_{\delta}(\alpha,\beta)\sim \Bigg\{\begin{array}{rcl} \delta^{1/2},&\mbox{for}& \alpha+\beta<1\\
	\delta^{1/2}(\log\frac{1}{\delta})^{1/2}, &\mbox{for}&\alpha+\beta=1\\
	\delta^{(2-\alpha-\beta)/2}, &\mbox{for}& \alpha+\beta>1.\end{array}\]
Next, we reduce the proof of Theorem~\ref{mainresult} to the following weighted $(2,2,1)$ estimate.  
\begin{theorem}\label{mainresult0}
	Let $n\geq3, \alpha_i\in [0,n-1)$ and $\beta_i\in [0,1)$ for $i=1,2$. Then, for $\lambda>\max\{\frac{\alpha_1+\beta_1}{2}+\frac{\alpha_2+\beta_2}{2}-1,\frac{\alpha_1+\beta_1}{2}-\frac{1}{2},\frac{\alpha_2+\beta_2}{2}-\frac{1}{2},0\}$, the following weighted $(2,2,1)$-estimate of the maximal bilinear cone multiplier operator $ \mathcal{T}^{\lambda}_{*}$  holds 
	\[\Vert \mathcal{T}^{\lambda}_{*}(f,g)\Vert_{L^{1}(\omega)}\lesssim \Vert f\Vert_{L^{2}(\omega_{\alpha_1,\beta_1})} \Vert g\Vert_{L^{2}(\omega_{\alpha_2,\beta_2})},\]
	where $\omega(x)=\omega^{1/2}_{\alpha_1,\beta_1}(x)\omega^{1/2}_{\alpha_2,\beta_2}(x)$.\\
	Consequently, for $f\in L^{2}(\omega_{\alpha_1,\beta_1}), g\in L^{2}(\omega_{\alpha_2,\beta_2})$ and the same range of parameters as above,  the pointwise almost everywhere convergence
	\[\lim_{R\to\infty}\mathcal{T}^{\lambda}_{R}(f,g)(x)=f(x)g(x)~~~\text{for a.e. }x\in\mathbb{R}^{n},\]
	holds, where \[supp(\hat{f}), supp(\hat{g})\subset \{\xi\in\mathbb{R}^n: \frac{1}{2}<|\xi_n|<2\}.\]
\end{theorem}
We build the proof of Theorem~\ref{mainresult0} with suitable modifications in the techniques developed by Jotsaroop and Shrivastava~\cite{JS} and Chen, He, Li and Yan~\cite{CHLY}. 

\subsection*{Organization of paper} 
Section \ref{3} is devoted to the proof of Theorem \ref{mainresult} and Section \ref{mainproof} contains the proof of Theorem \ref{mainresult0}. 
In the Appendix (Section \ref{Appendix}) we point out a direct proof of $L^2(\mathbb{R}^n)\times L^2(\mathbb{R}^n)\to L^1(\mathbb{R}^n)$ boundedness of the operator $\mathcal{T}^{\lambda}_{*}$.
\section{Theorem \ref{mainresult0} implies Theorem~\ref{mainresult}}\label{3}
The following result provides the key argument to complete this assertion that Theorem \ref{mainresult0} implies Theorem~\ref{mainresult}.
	\begin{proposition}[\cite{CHLY}]\label{chen}
		Let $p\geq2$, $\alpha_1=\beta_1=0,$ $0\leq\alpha_3<(n-1)(1-\frac{2}{p})<\alpha_2,\alpha_4<1+(n-1)(1-\frac{2}{p})$, and $0\leq \beta_2<1-\frac{2}{p}<\beta_3$, $\beta_4<2-\frac{2}{p}$. Then, any $f\in L^{p}(\mathbb{R}^n)$ can be decomposed as \[f=f_1+f_2+f_3+f_4\]
		such that $f_i\in L^{2}(\omega_{\alpha_i,\beta_i})$ for $i=1,2,3,4$. Moreover, if $supp(\hat{f})\subset \{\xi: \frac{1}{2}<|\xi_n|<2\}$, then $supp(\hat{f_i})\subset \{\xi: \frac{1}{10}<|\xi_n|<10\}$.
	\end{proposition}
	Let $\lambda>0$ and $2\leq p_1,p_2<\frac{2n}{n-1-2\lambda}$ be such that $1/p=1/p_1+1/p_2$. Invoke~Proposition~\ref{chen} to decompose the functions as 
	\[f=\sum^{4}_{i=1}f_{i},~g=\sum^{4}_{i=1}g_i,\] 
	where $f_i\in L^{2}(\omega_{\alpha_1^i,\beta_1^i})$ and $g_i\in L^{2}(\omega_{\alpha_2^i,\beta_2^i})$ with 
	\[\alpha^{i}_j+\beta^i_j<n(1-\frac{2}{p_j})+\epsilon,~\text{for}~j=1,2~\text{and}~ i=1,2,3,4. \] 
	Since $\lambda>\max\{\frac{\alpha^i_1+\beta^i_1}{2}+\frac{\alpha^i_2+\beta^i_2}{2}-1,\frac{\alpha_1+\beta_1}{2}-\frac{1}{2},\frac{\alpha_2+\beta_2}{2}-\frac{1}{2},0\}$ and $\epsilon>0$ is arbitrarily small, applying Theorem \ref{mainresult0} to each pair $(f_i, g_j)$, we deduce the pointwise almost everywhere convergence 
	\[\lim_{R\to\infty}\mathcal{T}^{\lambda}_{R}(f,g)(x)=f(x)g(x)~~~\text{for a.e. }x\in\mathbb{R}^{n},\]
	for 
	\[\lambda>\max\Big\{n(1-\frac{1}{p})-1,n(\frac{1}{2}-\frac{1}{p_1})-\frac{1}{2}, n(\frac{1}{2}-\frac{1}{p_2})-\frac{1}{2},0\Big\}.\]
	This proves Theorem~\ref{mainresult} under the assumption that Theorem~\ref{mainresult0} holds. \qed

\section{Proof of Theorem \ref{mainresult0}}\label{mainproof}
Recall that the bilinear cone multiplier is given by  $$m^{\lambda}(\frac{\xi'}{R\xi_n},\frac{\eta'}{R\eta_n})=\Big(1-\frac{|\xi'|^2}{R^2\xi^2_n}-\frac{|\eta'|^2}{R^2\eta^2_n}\Big)^{\lambda}_{+}\varphi(\xi_n)\varphi(\eta_n), \lambda>0.$$ 

Let $\psi\in C^{\infty}_{c}([1/2,2])$ and $\psi_{1}\in C^{\infty}_{c}([-3/4,3/4])$ be such that 
\[\sum_{j\geq2}\psi(2^{j}(1-t))+\psi_{1}(t)=1,~ t\in [0,1].\]
We begin with the following decomposition of the bilinear multiplier $m^{\lambda}$ with respect to the partition of unity as above in the $\xi$ variable. Consider  
\begin{eqnarray*}
	m^{\lambda}(\frac{\xi'}{R\xi_n},\frac{\eta'}{R\eta_n})&=&\sum_{j\geq2}\psi\Big(2^{j}(1-\frac{|\xi'|^2}{R^2\xi^2_n})\Big)m^{\lambda}(\frac{\xi'}{R\xi_n},\frac{\eta'}{R\eta_n})+\psi_{1}(\frac{|\xi'|^2}{R^2\xi^2_n})m^{\lambda}(\frac{\xi'}{R\xi_n},\frac{\eta'}{R\eta_n})\\
	&=&\sum_{j\geq2}m^{\lambda,j}(\frac{\xi'}{R\xi_n},\frac{\eta'}{R\eta_n})+m^{\lambda,1}(\frac{\xi'}{R\xi_n},\frac{\eta'}{R\eta_n}),
\end{eqnarray*}
where \[m^{\lambda,j}(\frac{\xi'}{R\xi_n},\frac{\eta'}{R\eta_n})=\psi\Big(2^{j}(1-\frac{|\xi'|^2}{R^2\xi^2_n})\Big)m^{\lambda}(\frac{\xi'}{R\xi_n},\frac{\eta'}{R\eta_n}),~j\geq 2~~\text{and}~\]
\[m^{\lambda,1}(\frac{\xi'}{R\xi_n},\frac{\eta'}{R\eta_n})=\psi_{1}(\frac{|\xi'|^2}{R^2\xi^2_n})m^{\lambda}(\frac{\xi'}{R\xi_n},\frac{\eta'}{R\eta_n}).\]
Let $\mathcal{T}^{\lambda}_{R,j}$ denote the bilinear multiplier operator associated with $m^{\lambda,j}$, that is, it is given by  \[\mathcal{T}^{\lambda}_{R,j}(f,g)(x):=\int_{\mathbb{R}^{2n}}m^{\lambda,j}\left(\frac{\xi'}{R\xi_n},\frac{\eta'}{R\eta_n}\right) \hat{f}(\xi)\hat{g}(\eta)e^{2\pi\iota x\cdot(\xi+\eta)}~d\xi d\eta,~~ j\geq 1.\]
This gives us the following decomposition of the maximal bilinear cone multiplier operator
\begin{align} \label{maindecom}
	\mathcal{T}^{\lambda}_{*}(f,g)(x)
	&\leq \sum_{j=1}^{\infty} \mathcal{T}^{\lambda}_{*,j}(f,g)(x),
	\end{align}
	where 
\begin{align*} \mathcal{T}^{\lambda}_{*,j}(f,g)(x)&=\sup_{R>0}|\mathcal{T}^{\lambda}_{R,j}(f,g)(x)|. 
\end{align*}

We prove suitable weighted $(2,2,1)$ estimates of the maximal operator $\mathcal{T}^{\lambda}_{*,j}$ for each $j\geq 1.$ We begin with the terms with $j\geq 2$. The remaining case of $j=1$ will require another decomposition using the partition of unity as earlier with respect to the other variable $\eta$. 

Therefore, let us first consider the case of $j\geq 2.$ We make use of the following  identity, see Stein and Weiss from [page 278,\cite{SW}], to decompose the multiplier $m^{\lambda,j}$ further. 

\begin{align}\label{SteinWeiss} \Big(1-\frac{|m|^2}{R^2}\Big)^{\lambda}=c_{\lambda} R^{-2\lambda}\int^{R}_{|m|}(R^2-t^2)^{\mu-1} t^{2\nu+1} \Big(1-\frac{|m|^2}{t^2}\Big)^{\nu}~dt,
	\end{align}
where $\lambda>0$ and $\lambda=\mu+\nu$ with $\mu>0,\nu>-1$ and the constant $c_{\lambda}=\frac{2\Gamma(\mu+\nu+1)}{\Gamma(\nu+1)\Gamma(\mu)}$. 
The identity above yields 

\begin{eqnarray*}
	&&m^{\lambda,j}(\frac{\xi'}{R\xi_n},\frac{\eta'}{R\eta_n})\\
	&&=c_{\lambda}R^{-2\lambda}\varphi(\xi_n)\varphi(\eta_n)\psi(2^{j}(1-\frac{|\xi'|^2}{R^2\xi^2_n}))\int^{R\sqrt{2^{1-j}}}_{0}(R^2\phi_{R}(\xi)-t^2)^{\mu-1}_{+} t^{2\nu+1}(1-\frac{|\eta'|^2}{t^2\eta^2_n})^{\nu}_{+}~dt,
\end{eqnarray*}
where $\lambda=\mu+\nu$ with $\mu>0$, $\nu>-1$ and $\phi_{R}(\xi)=(1-\frac{|\xi'|^2}{R^2\xi^2_n})_{+}$.

Define
\begin{eqnarray*}
	B^{R,t}_{j,\mu}f(x)&=&\int_{\mathbb{R}^n}\psi\Big(2^{j}(1-\frac{|\xi'|^2}{R^2\xi^2_n})\Big)\varphi(\xi_n)(R^2-\frac{|\xi'|^2}{\xi^2_n}-t^2)^{\mu-1}_{+}\hat{f}(\xi) e^{2\pi\iota x\cdot\xi} d\xi,\\
	T^{\nu}_{t}g(x)&=&\int_{\mathbb{R}^n}(1-\frac{|\eta'|^2}{t^2\eta^2_n})^{\nu}_{+}\varphi(\eta_n) \hat{g}(\eta) e^{2\pi\iota x\cdot\eta} d\eta.
\end{eqnarray*}
Observe that for $j\geq 2$, the operator  $\mathcal{T}^{\lambda}_{R,j}$ can be expressed as 
\begin{eqnarray}\label{splitoperator}
	\mathcal{T}^{\lambda}_{R,j}(f,g)(x):= c_{\lambda}R^{-2\lambda}\int^{R\sqrt{2^{1-j}}}_{0}B^{R,t}_{j,\mu}f(x) T^{\nu}_{t}g(x)t^{2\nu+1}~dt.
\end{eqnarray}
We apply Cauchy--Schwarz inequality to get that  
\begin{eqnarray*}
	|\mathcal{T}^{\lambda}_{R,j}(f,g)(x)|\leq c_{\lambda}R^{-2\lambda}\Big(\int^{R\sqrt{2^{1-j}}}_{0}|B^{R,t}_{j,\mu}f(x)t^{2\nu+1}|^2~dt\Big)^{1/2} \Big(\int^{R\sqrt{2^{1-j}}}_{0}|T^{\nu}_{t}g(x)|^{2}~dt\Big)^{1/2}.
\end{eqnarray*}
Further, the change of variable $t\to Rt$ in the first expression of the inequality above yields that 
\begin{align}\label{operatorT_j}
	|\mathcal{T}^{\lambda}_{R,j}(f,g)(x)|&\leq c_{\lambda}2^{-j/4}\Big(\int^{\sqrt{2^{1-j}}}_{0}|S^{R,t}_{j,\mu}f(x)t^{2\nu+1}|^2~dt\Big)^{1/2} \Big(\frac{1}{R\sqrt{2^{1-j}}}\int^{R\sqrt{2^{1-j}}}_{0}|T^{\nu}_{t}g(x)|^{2}~dt\Big)^{1/2},
\end{align}
where 
\begin{eqnarray*}
	S^{R,t}_{j,\mu}f(x)=\int_{\mathbb{R}^n}\psi\Big(2^{j}(1-\frac{|\xi'|^2}{R^2\xi^2_n})\Big)\varphi(\xi_n)(1-\frac{|\xi'|^2}{R^2\xi^2_n}-t^2)^{\mu-1}_{+}\hat{f}(\xi) e^{2\pi\iota x\cdot\xi} ~d\xi.
\end{eqnarray*}
In the expression above we have adjusted the factor $\sqrt{2^{1-j}}$ with second integral. 

In view of the decomposition of the bilinear operator $\mathcal{T}^{\lambda}_{R,j}$ as above, the desired weighted $(2,2,1)$ estimate can be obtained by proving appropriate weighted $L^2$-estimates of the following maximal operators
\[f\rightarrow \sup_{R>0}\Big(\int^{\sqrt{2^{1-j}}}_{0}|S^{R,t}_{j,\mu}f(x)t^{2\nu+1}|^2~dt\Big)^{\frac{1}{2}}~~\text{and}~~g\rightarrow \sup_{R>0}\Big(\frac{1}{R\sqrt{2^{1-j}}}\int^{R\sqrt{2^{1-j}}}_{0}|T^{\nu}_{t}g(x)|^{2}~dt\Big)^{\frac{1}{2}}.\]
More precisely, we prove the following key estimates of these operators.
\begin{lemma}\label{SRtjbeta}
	Let $n\geq3$ and $\alpha\in [0,n-1), \beta\in [0,1)$. Then for $\mu>\max\{\frac{\alpha+\beta}{2},\frac{1}{2}\}$ we have 
	\[\Big\Vert \sup_{R>0}\Big(\int^{\sqrt{2^{1-j}}}_{0}|S^{R,t}_{j,\mu}f(x)t^{2\nu+1}|^2~dt\Big)^{1/2} \Big\Vert_{L^{2}(\omega_{\alpha,\beta})}\leq  C 2^{-j(\mu+\nu)} 2^{\frac{3j}{4}} A_{2^{-j}(\alpha,\beta)}\Vert f\Vert_{L^{2}(\omega_{\alpha,\beta})}.\]
\end{lemma}
\begin{lemma}\label{Tdeltaweightedestimate01}
	Let $n\geq3$ and $\alpha\in [0,n-1), \beta\in [0,1)$. Then for $\nu>\max\{\frac{\alpha+\beta}{2}-1,-\frac{1}{2}\}$
	\[\Big\Vert \sup_{R>0}\Big(\frac{1}{R\sqrt{2^{1-j}}}\int^{R\sqrt{2^{1-j}}}_{0}|T^{\nu}_{t}g(x)|^{2}~dt\Big)^{1/2}\Big\Vert_{L^{2}(\omega_{\alpha,\beta})}\leq C  \Vert g\Vert_{L^{2}(\omega_{\alpha,\beta})}.\] 
\end{lemma}
The use of H\"{o}lder's inequality along with results in Lemma \ref{SRtjbeta} and Lemma \ref{Tdeltaweightedestimate01} yields the following weighted $(2,2,1)$-estimate of the  maximal operator $\mathcal{T}^{\lambda}_{*,j}, j\geq 2.$
\begin{theorem} For $j\geq2$, the following estimate holds
\begin{eqnarray}\label{Tjestimate}
	\Vert \mathcal{T}^{\lambda}_{*,j}(f,g)\Vert_{L^{1}(\omega)}\lesssim 2^{-j/4} 2^{-j(\mu+\nu)}2^{\frac{3j}{4}} A_{2^{-j}}(\alpha_1,\beta_1) \Vert f\Vert_{L^{2}(\omega_{\alpha_1,\beta_1})} \Vert g\Vert_{L^{2}(\omega_{\alpha_2,\beta_2})},
\end{eqnarray}
where $\omega(x)=\omega^{1/2}_{\alpha_1,\beta_1}(x)\omega^{1/2}_{\alpha_2,\beta_2}(x)$, and $(\alpha_1,\beta_1)$, $(\alpha_2,\beta_2)$ satisfy the hypothesis of Lemma \ref{SRtjbeta} and Lemma \ref{Tdeltaweightedestimate01} respectively. 
\end{theorem}
As a consequence, we get the desired weighted $(2,2,1)$ estimate for the operator $\sum_{j\geq2}\mathcal{T}^{\lambda}_{*,j}$ for all  $\lambda>\max\{\frac{\alpha_1+\beta_1}{2}+\frac{\alpha_2+\beta_2}{2}-1,\frac{\alpha_1+\beta_1}{2}-\frac{1}{2},\frac{\alpha_2+\beta_2}{2}-\frac{1}{2},0\}$. 

Therefore, we need to prove Lemma \ref{SRtjbeta} and Lemma \ref{Tdeltaweightedestimate01}.
\subsection{Proof of Lemma~\ref{SRtjbeta}}\label{case2}
We require some new notation in order to prove this lemma. For $\Psi\in C_{c}^{\infty}([\frac{1}{4},1])$ and $0<\delta< 1/16$, consider the smooth cut-off of the linear cone multiplier operator 
	\[T^{\Psi}_{\delta,t}f(x)=\int_{\mathbb{R}^{n}}\varphi(\xi_n)\Psi\Big(\delta^{-1}(1-\frac{|\xi'|^2}{t^2\xi^{2}_n})\Big) \hat{f}(\xi) e^{2\pi\iota x\cdot\xi}~d\xi.\]
	Observe that the multiplier $\frak{m}(\xi')=\Psi(\delta^{-1}(1-|\xi'|^2))$ satisfies the following property  \[|\partial_{\xi'}^{L} \frak{m}(\xi')|\lesssim \delta^{-|L|},\]
	where the multi-index $L=(l_1,l_2,\dots,l_{n-1})$ with $l_{i}\in \mathbb{N}\cup \{0\}$. 
It is natural to decompose the cone multiplier into these smooth multipliers. We refer the reader to~\cite{CHLY} for more  details. The square function associated with $T^{\Psi}_{\delta,t}f$ is considered in~\cite{CHLY} and is defined by
\begin{align*}
	G_{\delta}(f)(x)=\Big(\int^{\infty}_{0}|T^{\Psi}_{\delta,t}f(x)|^2 \frac{dt}{t}\Big)^{1/2}.
\end{align*}
Invoking Proposition $3.3$ and Theorem $5.1$ from \cite{CHLY}, we have the following weighted $L^2$-estimate of the square function.
\begin{lemma}\label{Tpsiweightedestimate1}\cite{CHLY}
	Let $n\geq3$, $0<\delta\ll 1$  and $\alpha\in [0,n-1), \beta\in [0,1)$. Then 
	\[\Big\Vert \Big(\int^{\infty}_{0}|T^{\Psi}_{\delta,t}f(\cdot)|^2 \frac{dt}{t}\Big)^{1/2}\Big\Vert_{L^{2}(\omega_{\alpha,\beta})}\lesssim A_{\delta}(\alpha,\beta) \Vert f \Vert_{L^{2}(\omega_{\alpha,\beta})}.\]
\end{lemma}
The square function estimate in the lemma above allows us to deduce weighted $L^2$-estimates of the maximal function  $\sup_{t>0}|T^{\Psi}_{\delta,t}f|.$ More precisely, using Sobolev embedding theorem we can write 
	\[\sup_{t>0}|T^{\Psi}_{\delta,t}f(x)|^{2}\leq \delta^{-1} \Big(\int^{\infty}_{0}|T^{\Psi}_{\delta,t}f(x)|^2 \frac{dt}{t}\Big)^{1/2} \Big(\int^{\infty}_{0}|T^{\tilde{\Psi}}_{\delta,t}f(x)|^2 \frac{dt}{t}\Big)^{1/2},\]
	where $\tilde{\Psi}(\xi',\xi_n)=\delta \xi' \cdot \nabla_{n-1}\Psi(\xi',\xi_n) $. 
	
	Therefore, an application of Lemma \ref{Tpsiweightedestimate1} yields the following result for the maximal function.
\begin{lemma}\label{Tpsiweightedestimate2}
	Let $n\geq3$, $0<\delta\ll 1$  and $\alpha\in [0,n-1), \beta\in [0,1)$. Then 
	\[\Vert \sup_{t>0}|T^{\Psi}_{\delta,t}f(\cdot)|\Vert_{L^{2}(\omega_{\alpha,\beta})}\lesssim \delta^{-1/2}A_{\delta}(\alpha,\beta) \Vert f \Vert_{L^{2}(\omega_{\alpha,\beta})}.\]
\end{lemma}
With these auxiliary results, we are in a position to prove Lemma \ref{SRtjbeta}. Recall the definition of the operator $S^{R,t}_{j,\mu}$. 
	\begin{eqnarray*}
		S^{R,t}_{j,\mu}f(x)=\int_{\mathbb{R}^n}\psi\Big(2^{j}(1-\frac{|\xi'|^2}{R^2\xi^2_n})\Big)\varphi(\xi_n)(1-\frac{|\xi'|^2}{R^2\xi^2_n}-t^2)^{\mu-1}_{+}\hat{f}(\xi) e^{2\pi\iota x\cdot\xi} ~d\xi.
	\end{eqnarray*}
	Observe that the support of $\psi$ is in $[1/2,2]$. Consequently, when $t^{2}\in [0,2^{-j-1})$ or  $t<\sqrt{2^{-1-j}}$, the multiplier $(1-\frac{|\xi'|^2}{R^2\xi^2_n}-t^2)^{\mu-1}$ does not have singularity. Therefore, we consider the cases of $t<\sqrt{2^{-1-j}}$ and $t\geq \sqrt{2^{-1-j}}$ separately. 
	\subsection*{Case I:} Consider the case when $t\in [0,\sqrt{2^{-1-j-\epsilon_{0}}})$ for some small $\epsilon_{0}>0$. Write $\mu-1=\gamma$. If $-1<\gamma<0$, we choose $\gamma=-\rho$ for $\rho>0$. Using Taylor series expansion we write 
	\begin{eqnarray*}(1-\frac{|\xi'|^{2}}{R^2\xi^{2}_n}-t^{2})^{-\rho}&=&(1-\frac{|\xi'|^{2}}{R^2\xi^{2}_{n}})^{-\rho} \Big(1-\frac{t^2}{1-\frac{|\xi'|^{2}}{R^2\xi^{2}_{n}}}\Big)^{-\rho}\\
		&=&2^{j\rho}(2^{j}(1-\frac{|\xi'|^{2}}{R^2\xi^{2}_{n}}))^{-\rho}\sum_{k\geq0}\frac{\Gamma(\rho+k)}{\Gamma(\rho) k!} \Big( \frac{2^{j}t^2}{2^{j}(1-\frac{|\xi'|^{2}}{R^2 \xi^{2}_{n}})}\Big)^{k}.
	\end{eqnarray*}
	Observe that $\frac{t^{2}}{(1-\frac{|\xi'|^2}{R^2 \xi^{2}_{n}})}<2^{-j-1-\epsilon_{0}} 2^{j+1}<1$. Therefore, the power series as above converges and the operator $S^{R,t}_{j,\mu}$ takes the form 
	\begin{align}\label{Appen1}
		S^{R,t}_{j,\gamma+1}f(x)=2^{j\rho}\sum_{k\geq0} \frac{\Gamma(\rho+k)}{\Gamma(\rho)k!} (2^{j}t^{2})^{k} T_{2^{-j},R}^{\psi^{k}}f(x),
		\end{align}
	where $\psi^{k}(x)=x^{-k-\rho} \psi(x)$ belongs to the space $ C^{\infty}_{c}([1/2,2])$. 
	
	Invoking Lemma \ref{Tpsiweightedestimate2} we get that  
	\[\Vert \sup_{R>0}|S^{R,t}_{j,\gamma+1}f|\Vert_{L^{2}(\omega_{\alpha,\beta})}\lesssim 2^{j\rho} 2^{j/2} A_{2^{-j}}(\alpha,\beta)\sum_{k\geq0} \frac{\Gamma(\rho+k)}{\Gamma(\rho)k!} 2^{-(1+\epsilon_{0})k} k^{N} 2^{N+k} \Vert f\Vert_{L^{2}(\omega_{\alpha,\beta})}.\]
	Hence we get
	\begin{eqnarray*}
		\Big\Vert \sup_{R>0}\Big(\int^{\sqrt{2^{1-j}}}_{0}|S^{R,t}_{j,\mu}f(x)t^{2\nu+1}|^2~dt\Big)^{1/2} \Big\Vert_{L^{2}(\omega_{\alpha,\beta})}&\lesssim& 2^{j\rho} 2^{j/2} A_{2^{-j}}(\alpha,\beta) 2^{-j(4\nu+3)/4} \Vert f\Vert_{L^{2}(\omega_{\alpha,\beta})}\\
        &=& 2^{j(1-\mu)} 2^{j/2} A_{2^{-j}}(\alpha,\beta) 2^{-j(4\nu+3)/4} \Vert f\Vert_{L^{2}(\omega_{\alpha,\beta})}\\
        &=& 2^{-j(\mu+\nu)} 2^{3j/4} A_{2^{-j}}(\alpha,\beta)  \Vert f\Vert_{L^{2}(\omega_{\alpha,\beta})}.
	\end{eqnarray*}
	Next, we consider the case $\gamma\geq0$. Indeed, observe that, it is enough to prove the result for $\gamma\in [0,1].$ The remaining cases of $\gamma > 1$ can be dealt with similar arguments as those used for $\gamma\in [0,1].$
	
	Therefore, we consider the case $\gamma=0$ (i.e. $\mu=1$)  first. In this case, observe that $\sup_{R>0}|S^{R,t}_{j,1}(f)|$ is same as $\sup_{R>0}|T^{\Psi}_{\delta,R}(f)|$ with $\delta=2^{-j}$. Hence, we get
    \begin{eqnarray*}\Big\Vert \Big(\sup_{R>0}\int^{\sqrt{2^{-j-1-\epsilon_0}}}_{0}(S^{R,t}_{j,1}f(\cdot)t^{2\nu+1})^2 dt\Big)^{1/2}\Big\Vert_{L^{2}(\omega_{\alpha,\beta})}&\lesssim &2^{-j/4}2^{-j(\nu+\frac{1}{2})} 2^{j/2} A_{2^{-j}}(\alpha,\beta) \Vert f\Vert_{L^{2}(\omega_{\alpha,\beta})}\\
    &=& 2^{-j(1+\nu)} 2^{3j/4} A_{2^{-j}}(\alpha,\beta) \Vert f\Vert_{L^{2}(\omega_{\alpha,\beta})}.
    \end{eqnarray*}
    Next, for the case $\gamma=1$, we rewrite the underlying multiplier as 
    \begin{align*}
        \psi(2^{j}(1-\frac{|\xi'|^2}{R^2\xi^2_n}))(1-\frac{|\xi'|^2}{R^2\xi^2_n}-t^2)
        &&=2^{-j} \Big(2^{j}(1-\frac{|\xi'|^2}{R^2\xi^2_n})\Big)\psi(2^{j}(1-\frac{|\xi'|^2}{R^2\xi^2_n}))-t^2 \psi(2^{j}(1-\frac{|\xi'|^2}{R^2\xi^2_n})).
    \end{align*}
    Since $\Big(2^{j}(1-\frac{|\xi'|^2}{R^2\xi^2_n})\Big)\psi(2^{j}(1-\frac{|\xi'|^2}{R^2\xi^2_n}))$ behaves like $\psi(2^{j}(1-\frac{|\xi'|^2}{R^2\xi^2_n}))$, we get that  
    \[|S^{R,t}_{j,2}f(x)|\leq 2^{-j} |T^{\Psi}_{2^{-j},R}f(x)|+t^2 |T^{\Psi}_{2^{-j},R}f(x)|.\] 
    This gives us the desired estimate. 
    
    Finally, for $0<\gamma<1$, we write $\gamma=1+\zeta$ with $\zeta\in (-1,0)$. Then 
    \begin{eqnarray*}
        &&\psi(2^{j}(1-\frac{|\xi'|^2}{R^2\xi^2_n}))(1-\frac{|\xi'|^2}{R^2\xi^2_n}-t^2)^{\gamma}\\
        &&= 2^{-j}\Psi_1(\frac{|\xi'|^2}{R^2\xi^2_n})(1-\frac{|\xi'|^2}{R^2\xi^2_n}-t^2)^{\zeta}-t^2 \Psi_2(\frac{|\xi'|^2}{R^2\xi^2_n})(1-\frac{|\xi'|^2}{R^2\xi^2_n}-t^2)^{\zeta},
    \end{eqnarray*}
    where $\Psi_1(\frac{|\xi'|^2}{R^2\xi^2_n})=2^j (1-\frac{|\xi'|^2}{R^2\xi^2_n}) \psi(2^j(1-\frac{|\xi'|^2}{R^2\xi^2_n}))$  and $\Psi_2(\frac{|\xi'|^2}{R^2\xi^2_n})=\psi(2^j(1-\frac{|\xi'|^2}{R^2\xi^2_n}))$. Observe that the boundedness of the square function associated with $\Psi_2(\frac{|\xi'|^2}{R^2\xi^2_n})(1-\frac{|\xi'|^2}{R^2\xi^2_n}-t^2)^{\zeta}$ can be handled similarly as the case of $-1<\gamma<0$. On the other hand, note that the function $\Psi_1(\frac{|\xi'|^2}{R^2\xi^2_n})$ behaves like  $\psi(2^j(1-\frac{|\xi'|^2}{R^2\xi^2_n}))$. Therefore, the square function associated with $\Psi_1(\frac{|\xi'|^2}{R^2\xi^2_n})$ can also be dealt with similarly as in the case of $-1<\gamma<0$.  This completes the proof of desired bounds for the ${\bf Case~ I}$.
     
	\subsection*{Case II:} Let $t^2\in [2^{-j-1-\epsilon_{0}}, 2^{-j+1}]$. Consider 
	\begin{eqnarray*}
		\psi\Big(2^{j} (1-\frac{|\xi'|^2}{R^2\xi^2_n})\Big)(1-\frac{|\xi'|^2}{R^2\xi^2_n}-t^2)^{\mu-1}_{+}=(1-t^{2})^{\mu-1} \psi\Big(2^{j} (1-\frac{|\xi'|^2}{R^2\xi^2_n})\Big)\Big(1-\frac{|\xi'|^2}{(1-t^{2})\xi^{2}_nR^2}\Big)^{\mu-1}_{+}.
	\end{eqnarray*}
	Applying the change of variable  $1-t^2 =s^2$, we get that  
	\[\sup_{R>0}\Big(\int^{\sqrt{2^{1-j}}}_{\sqrt{2^{-j-1-\epsilon_0}}}|S^{R,t}_{j,\mu}f(x)|^{2}~tdt\Big)^{1/2}\simeq \sup_{R>0}\Big(\int^{\sqrt{1-2^{-j-1-\epsilon_{0}}}}_{\sqrt{1-2^{-j+1}}}|S'^{R,s}_{j,\mu}f(x)|^{2}~ds\Big)^{1/2},\]
	where 
	\[S'^{R,s}_{j,\mu}f(x)=\int \varphi(\xi_n)\psi\Big(2^{j}(1-\frac{|\xi'|^2}{R^{2}\xi^{2}_{n}})\Big) (1-\frac{|\xi'|^2}{s^2 R^2 \xi^{2}_{n}})^{\mu-1}_{+} \hat{f}(\xi) e^{2\pi\iota x\cdot\xi} d\xi.\]
    Hence,  
    \begin{align}\label{Appen2}\sup_{R>0}\Big(\int^{\sqrt{2^{1-j}}}_{\sqrt{2^{-j-1-\epsilon_0}}}|S^{R,t}_{j,\mu}f(x)|^{2}~ t^{4\nu+1} tdt\Big)^{1/2}&\simeq 
 2^{-j\nu} 2^{-j/4}\sup_{R>0}\Big(\int^{\sqrt{1-2^{-j-1-\epsilon_{0}}}}_{\sqrt{1-2^{-j+1}}}|S'^{R,s}_{j,\mu}f(x)|^{2}~ds\Big)^{1/2}.
 \end{align}
	The following result holds for the operator $S'^{R,s}_{j,\mu}$ when $2\leq j<M'$. 
	\begin{proposition}\label{prop5.3}
		Let  $2\leq j<M', \alpha\in [0,n-1), \beta\in [0,1)$. Then for $\mu>\max\{\frac{\alpha+\beta}{2},\frac{1}{2}\}$, we have 
		\[\Big\Vert \sup_{R>0}\Big(\int^{s_2}_{s_1}|S'^{R,s}_{j,\mu}f(x)|^2~ds\Big)^{1/2}\Big\Vert_{L^{2}(\omega_{\alpha,\beta})}\leq C_{M',n} \Vert f\Vert_{L^{2}(\omega_{\alpha,\beta})}.\]
	\end{proposition}

	\begin{proof}
		Note that we have 
		\begin{align*}S'^{R,s}_{j,\mu}f(x)&=T^{\psi}_{2^{-j},R}\circ T^{\mu-1}_{sR}f(x),
			\end{align*}
		where \[T^{\psi}_{2^{-j},R}f(x)=\int \varphi(\xi_n)\psi\Big(2^{j}(1-\frac{|\xi'|^2}{R^{2}\xi^{2}_{n}})\Big)\hat{f}(\xi) e^{2\pi\iota x\cdot\xi} d\xi.\]
Since we are dealing with finitely many values of $j$ here, we make use of the following crude estimate
        \begin{eqnarray}\label{Maximalfunction1}
        \sup_{R>0}|T^{\psi}_{2^{-j},R}f(x)|\lesssim C_j \mathfrak{M}f(x),~\text{a.e.}~x,
        \end{eqnarray}
        where \[\mathfrak{M}f(x',x_n):=\sup_{\mathfrak{Q}\ni (x',x_n)} \frac{1}{|\mathfrak{Q}|}\int_{\mathfrak{Q}} |f(y)|~dy,\] 
        and $\mathfrak{Q}=Q\times I$, $Q$ is cube in $\mathbb{R}^{n-1}$ with sides parallel to the coordinate axes, and $I$ is an interval in $\mathbb{R}$.
        Indeed, consider $K_{R,j}$ be the kernel of $T^{\psi}_{2^{-j},R}$. Then
        \begin{eqnarray*}
            K_{R,j}(y',y_n)&=&\int_{\mathbb{R}^n} \varphi(\xi_n)\psi\Big(2^{j}(1-\frac{|\xi'|^2}{R^{2}\xi^{2}_{n}})\Big) e^{2\pi\iota (y'\cdot \xi'+y_n\xi_n)} d\xi\\
            &=&\int_{\mathbb{R}} \varphi(\xi_n)e^{2\pi\iota y_n\xi_n} \int_{\mathbb{R}^{n-1}}\psi\Big(2^{j}(1-\frac{|\xi'|^2}{R^{2}\xi^{2}_{n}})\Big) e^{2\pi\iota y'\cdot \xi'} d\xi' d\xi_n\\
            &=&R^{n-1} \int_{\mathbb{R}}\xi^{n-1}_n \varphi(\xi_n)e^{2\pi\iota y_n\xi_n} \int_{\mathbb{R}^{n-1}}\psi\Big(2^{j}(1-|\xi'|^2)\Big) e^{2\pi\iota R\xi_n y'\cdot \xi'} d\xi' d\xi_n.
        \end{eqnarray*}
        
        Since $\psi$ and $\varphi$ are compactly supported smooth functions, we can easily verify the estimate of the kernel $K_{R,j}(y',y_n)$ to show that the the equation  \eqref{Maximalfunction1} holds with a constant $C_j$. Even though $C_j$ depends on $j$, it is sufficient for our purpose in this case.

It is well known that $\mathfrak{M}$ is bounded in $L^{p}(|x'|^{-\alpha}|x_n|^{-\beta})$ for $1<p<\infty$ with $1-n<-\alpha<(n-1)(p-1), -1<-\beta<p-1$.
   
		For a large number $n_0$ such that $\mu-1+n_0>n/2$, we have 
		\begin{eqnarray*}
			&&\sup_{R>0}\Big(\int^{s_2}_{s_1}|S'^{R,s}_{j,\mu}f(x)|^2~ds\Big)^{1/2}\\
			&&\leq \sum^{n_0}_{k=1}\sup_{R>0}\Big(\int^{s_2}_{s_1} |T^{\psi}_{2^{-j},R}\circ( T^{\mu-1+k}_{sR}-T^{\mu+k-2}_{sR})f(x)|^2\Big)^{1/2}+\sup_{R>0}\Big(\int^{s_2}_{s_1} |T^{\psi}_{2^{-j},R}\circ T^{\mu-1+n_0}_{sR}f(x)|^2\Big)^{1/2}\\
			&&\lesssim 2^{M'n}\Big(\int^{\infty}_{0}|\mathfrak{M}\circ \tilde{T}^{\mu-1}_{s}f(x)|^2~\frac{ds}{s}\Big)^{1/2}+ 2^{M'n} \mathfrak{M}\circ T^{\mu-1+n_0}_{*}f(x),
		\end{eqnarray*}
        where $\widehat{\tilde{T}^{\mu-1}_{s}f}(\xi)=\frac{|\xi'|^2}{s^2\xi^2_n}(1-\frac{|\xi'|^2}{s^2\xi^2_n})^{\mu-1}\hat{f}(\xi)$.
		Note that $ T^{\mu-1+n_0}_{*}f(x)\lesssim M_{HL}f(x)$ as $\mu-1+n_0>n/2$. 
		Therefore, invoking weighted boundedness of $\mathfrak{M}$ and $M_{HL}$, the Hardy--Littlewood maximal function 
		and Lemma \ref{Tdeltaweightedestimate01} we get the desired  $L^{2}(\omega_{\alpha,\beta})$ estimate of the operator under consideration  for $\mu>\max\{\frac{\alpha+\beta}{2},\frac{1}{2}\}$.
	\end{proof}
	
	For large values of $j$, we have the following estimate. 
	\begin{proposition}\label{prop5.4}
		Let $n\geq3$, $\mu>\frac{1}{2}$ and $j\geq M'$ for some large $M'\in\mathbb{N}$. Then for $\alpha\in [0,n-1), \beta\in [0,1)$ we have 
		\[\Big\Vert \sup_{R>0}\Big(\int^{s_2}_{s_1}|S'^{R,s}_{j,\mu}f(x)|^2~ds\Big)^{1/2} \Big\Vert_{L^{2}(\omega_{\alpha,\beta})}\lesssim 2^{-j(\mu-1)}  A_{2^{-j}}(\alpha,\beta)  \Vert f\Vert_{L^{2}(\omega_{\alpha,\beta})}.\]
	\end{proposition}
    Note that the above proposition implies 
    \[\Big\Vert\sup_{R>0}\Big(\int^{\sqrt{2^{1-j}}}_{\sqrt{2^{-j-1-\epsilon_0}}}|S^{R,t}_{j,\mu}f(x) t^{2\nu+1}|^{2}~dt\Big)^{1/2}\Big\Vert_{L^{2}(\omega_{\alpha,\beta})}\lesssim  2^{-j\nu} 2^{-j/4}2^{-j(\mu-1)}  A_{2^{-j}}(\alpha,\beta)  \Vert f\Vert_{L^{2}(\omega_{\alpha,\beta})}. \]
	\begin{proof}
		Consider a partition of unity to decompose the multiplier  as follows
		\begin{align*}
			\psi(2^j(1-\frac{|\xi'|^2}{R^2\xi^2_n}))(1-\frac{|\xi'|^2}{s^2R^2\xi^2_n})^{\gamma}_{+}&=&\sum_{k\geq2}2^{-k\gamma}\psi(2^j(1-\frac{|\xi'|^2}{R^2\xi^2_n}))\tilde{\psi}\Big(2^k(1-\frac{|\xi'|^2}{s^2R^2\xi^2_n})\Big)\\
			&&+\psi(2^j(1-\frac{|\xi'|^2}{R^2\xi^2_n}))\tilde{\psi_1}(\frac{|\xi'|^2}{s^2R^2\xi^2_n}),
		\end{align*}
		where $\tilde{\psi}\in C^{\infty}_c([1/2,2])$ and $\tilde{\psi_1}\in C^{\infty}_c([0,3/4])$.
		
		Observe that the support condition on $\psi$ and $\tilde{\psi}$ imply that $k\geq j-2$.

Next, we decompose the support of $\tilde{\psi}$ into annuli of small width so that $\psi(2^j(1-\frac{|\xi'|^2}{R^2\xi^2_n}))$ behaves like a constant when restricted to these annular regions. For, let $\phi\in C^{\infty}_c([-1,1])$ be such that $\sum_{l\in\mathbb{Z}}\phi(x-l)=1$ for $x\in \mathbb{R}$ and 
$$\sup_{\substack { x\in [-1,1]}}|\frac{d^l\phi(x)}{dx^l}|\leq L,~~0\leq l\leq N.$$ 

 Let $a=R^2 s^2(1-2^{-k+1})$ and $b=R^2 s^2(1-2^{-k-1})$, then for $\delta=\frac{2}{3}2^{-k\epsilon}$ with $0<\epsilon<1$, using the identity above on the interval $[a,b]$ we get that 
		\[\tilde{\psi}(2^k(1-\frac{|\xi'|^2}{R^2 s^2\xi^2_n}))=\sum_{0\leq l\leq [\delta^{-1}]+1}\phi\Big(2^{(1+\epsilon)k}(\frac{|\xi'|^2}{R^2 s^2\xi^2_n}-1+2^{-k+1}-2^{-(1+\epsilon)k}l)\Big)\tilde{\psi}(2^k(1-\frac{|\xi'|^2}{R^2 s^2\xi^2_n})).\]
		Since 
		\[\tilde{\psi}(2^k(1-\frac{|\xi'|^2}{R^2 s^2\xi^2_n}))=\int_{\mathbb{R}}\widehat{\tilde{\psi}}(t)e^{2\pi\iota t(2^k(1-\frac{|\xi'|^2}{R^2 s^2\xi^2_n}))}dt.\]
		We reorganize the exponential term in the expression above and use Taylor series expansion in the following way.  
		\[e^{2\pi\iota t(2^k(1-\frac{|\xi'|^2}{R^2 s^2\xi^2_n}))}=e^{2\pi\iota t2^k(2^{-k+1}-2^{-k(1+\epsilon)}l)} e^{-2\pi\iota t2^{-\epsilon k}2^{(1+\epsilon)k}(\frac{|\xi'|^2}{R^2s^2\xi^2_n}-1+2^{-k+1}-2^{-(1+\epsilon)k}l)}.\]
		We observe that the quantity $d^{\epsilon}_{k,l}=2^k(2^{-k+1}-2^{-k(1+\epsilon)}l)$ is uniformly bounded for $l\leq 1+[\delta^{-1}]$ and $k\geq2$.  Let  $\xi^{Rs}_{k,l}=(\frac{|\xi'|^2}{R^2s^2\xi^2_n}-1+2^{-k}d^{\epsilon}_{k,l})$ and write 
		\[e^{-2\pi\iota t2^{k}\xi^{Rs}_{k,l}}=\sum^{N-1}_{q=0}\frac{(-1)^q}{q!} (2\pi\iota 2^{-\epsilon k})^q 2^{q(1+\epsilon)k}(\xi^{Rs}_{k,l})^q t^q+r_{N}(2^k t\xi^s_{k,l}).\]
		Observe that \[\sup_{0\leq q\leq N-1}|\frac{d^q\tilde{\psi}}{dx^q}(d^{\epsilon}_{k,l})|\leq C,\] where C is independent of $k,l$. The discussion as above gives us the following form of the operator $S'^{R,s}_{j,\gamma+1}$
		\begin{align*}
			S'^{R,s}_{j,\gamma+1}f(x)&=\sum_{k\geq j-2} 2^{-k\gamma}\sum_{l\leq [\delta^{-1}]+1} \sum^{N-1}_{q=0} \frac{(-1)^q}{q!} (2\pi\iota 2^{-\epsilon k})^q \frac{d^q\tilde{\psi}}{dx^q}(d^{\epsilon}_{k,l}) U^{\phi_q,j}_{Rs, R}(k,l)f(x)\\
			&+\sum_{k\geq j-2}2^{-k\gamma} \sum_{l\leq [\delta^{-1}]+1} \int_{\mathbb{R}} \widehat{\tilde{\psi}} (t) e^{2\pi\iota d^{\epsilon}_{k,l}t} P^{j,k}_{Rs,R}(t,l)f(x) ~dt,
		\end{align*}
		where 
		\begin{align*}
			U^{\phi_q,j}_{Rs, R}(k,l)f(x)&=\int_{\mathbb{R}^n} \phi_{q}( 2^{(1+\epsilon)k}\xi^{Rs}_{k,l}) \psi\Big(2^{j}(1-\frac{|\xi'|^2}{R^2\xi^2_n})\Big) \hat{f}(\xi) e^{2\pi\iota x\cdot\xi} ~d\xi
		\end{align*}
		and 
	\begin{align*}
			P^{j,k}_{Rs,R}(t,l)f(x)&=\int_{\mathbb{R}^n} \phi (2^{(1+\epsilon)k}\xi^{Rs}_{k,l}) \psi\Big(2^{j}(1-\frac{|\xi'|^2}{R^2 \xi^2_n})\Big) r_{N}(2^k\xi^{Rs}_{k,l}t) \hat{f}(\xi) e^{2\pi\iota x\cdot \xi} ~d\xi.
		\end{align*}
		Therefore, the desired weighted $L^2$-estimate of the operator $S'^{R,s}_{j,\gamma+1}$ will be proved once we have suitable estimates of the operators $U^{\phi_q,j}_{Rs, R}(k,l)$ and $P^{j,k}_{Rs,R}(t,l)$. We discuss them below. 
		\subsection*{Estimates of the operator $U^{\phi_q,j}_{Rs, R}(k,l), 0\leq q\leq N-1$:} 
		It is enough to consider the case of $q=0$. The other cases  $1\leq q\leq N-1$ can be handled similarly. 
		As earlier, we decompose the multiplier of the operator $U^{\phi_0,j}_{Rs, R}(k,l)$ into multipliers supported on thin annular regions $\{\xi:\frac{|\xi'|^2}{R^2\xi^2_n}\in [a+s^22^{-(1+\epsilon)k}(l-1),a+s^22^{-(1+\epsilon)k}(l+1)]\},$ for $l\leq 1+[\delta^{-1}]$ and  $a=s^2(1-2^{-k+1})$. 
		Consider \[\psi(2^j(1-\frac{|\xi'|^2}{R^2\xi^2_n}))=\int_{\mathbb{R}}\hat{\psi}(u)e^{2\pi\iota(2^j(1-\frac{|\xi'|^2}{R^2\xi^2_n}))u} du.\] 
		We rewrite the exponential term as  
		\[e^{2\pi\iota(2^j(1-\frac{|\xi'|^2}{R^2\xi^2_n}))u}=e^{2\pi\iota 2^j(1-a-s^22^{-(1+\epsilon)k}l)u} e^{-2\pi\iota2^js^2(\frac{|\xi'|^2}{s^2R^2\xi^2_n}-1+2^{-k+1}-2^{-(1+\epsilon)k}l)u}\]
	and use Taylor series expansion of $e^{-2\pi\iota2^j s^2\xi^{Rs}_{k,l}u}$ to get that  
		\begin{align*}
			U^{\phi_0,j}_{Rs,R}(k,l)f(x)&=\sum^{N-1}_{m=0}\frac{1}{m!} s^2 2^{-d(1+\epsilon)m} 2^{-j\epsilon m}\psi^{m}(2^j(1-a-s^22^{-(1+\epsilon)k}l)) V^{\phi_{m}}_{Rs}(k,l)f(x)+X^{k}_{N,R,s}f(x), 
		\end{align*}
		where  $\phi_{m}(x)=x^m \phi(x)$, 
		\begin{align*}
			V^{\phi_m}_{Rs}(k,l)f(x)&=\int_{R^n} \phi_{m}(2^{(1+\epsilon)k}\xi^{Rs}_{k,l}) \hat{f}(\xi) e^{2\pi\iota x\cdot \xi} d\xi,\\
			X^{k}_{N,R,s}f(x)&=\int_{\mathbb{R}} \hat{\psi}(u) e^{2\pi\iota 2^j(1-a-s^22^{-(1+\epsilon)k}l)u} Q^{k}_{N,R,s}(u)f(x) du,
		\end{align*}
		with 
		\begin{align*}Q^{k}_{N,R,s}(u)f(x)=\int_{\mathbb{R}^n}\phi_{m} (2^{(1+\epsilon)k}\xi^{Rs}_{k,l}) r_{N}(2^j s^2\xi^{Rs}_{k,l}u) \hat{f}(\xi) e^{2\pi\iota x\cdot\xi} d\xi.
			\end{align*}
		Observe that change of variables argument yields that 
        \begin{eqnarray}\label{Vphi}
        \Big(\int^{s_2}_{s_1}|V^{\phi_{m}}_{Rs}(k,l)f(\cdot)|^2 ds\Big)^{1/2}
			&\lesssim \Big( \int^{\infty}_{0}|T^{\phi_{m}}_{\varrho,s}f(x)|^{2} \frac{ds}{s}\Big)^{1/2}.
			\end{eqnarray}
		Note that using Lemma \ref{Tpsiweightedestimate1} we get the following estimate
		\[\Big\Vert \Big(\int^{\infty}_{0}|T^{\phi_{m}}_{\varrho,s}f(x)|^{2} \frac{ds}{s}\Big)^{1/2}\Big\Vert_{L^{2}(\omega_{\alpha,\beta})}\lesssim  A_{\varrho}(\alpha,\beta)\Vert f\Vert_{L^{2}(\omega_{\alpha,\beta})},\]
		where $\varrho\simeq 2^{-(1+\epsilon)k}$. 
	\end{proof}
	Next, for the operator $X^{k}_{N,R,s}$, we have the following estimate. 
	\begin{eqnarray}\label{estX}
		\Big\Vert \sup_{R>0,s\in [s_1,s_2]}|X^{k}_{N,R,s}f(x)|\Big\Vert_{L^{2}(\omega_{\alpha,\beta})}\lesssim A_{2^{-k}}(\alpha,\beta)\Vert f\Vert_{L^{2}(\omega_{\alpha,\beta})},
	\end{eqnarray}
where $\alpha\in [0,n-1), \beta\in [0,1)$. 

Indeed, by a change of variables argument, the operator under consideration can be shown to be dominated by the operator $\mathfrak{M}$. Hence, the desired weighted $L^2$-estimate follows. 
\subsection{Estimates of the operator $P^{j,k}_{Rs,R}$:}\label{Pjk} In order to deal with the operator $P^{j,k}_{Rs,R}$,  we decompose the support of $\psi$ further  into small annuli of width $2^{-(1+\epsilon)m}$ such that 
\begin{align*}
    P^{j,k}_{Rs,R}(t,l)f(x)&=\sum^{N-1}_{m=0} \frac{1}{m!} s^2 2^{-d(1+\epsilon)m} 2^{-j\epsilon m} \psi_{j}(k,l,m) D^{m,k}_{s,R}(t,l)f(x)\\
    &+\int_{\R}\int_{\mathbb{R}^n} \hat{\psi}(t') s^{R}_{N}(t,t') \varphi(2^{(1+\epsilon)k}\xi^{Rs}_{k,l})\hat{f}(\xi)e^{2\pi\iota x\cdot \xi} e^{2\pi\iota 2^j t'(1-a-s^22^{-(1+\epsilon)k}l)}  d\xi dt',
\end{align*}
where 
$\psi_{j}(k,l,m)=\frac{d^m\psi}{dx^m}\Big(2^j(1-a-s^2 2^{-(1+\epsilon)k}l)\Big), 
    ~~s^{R}_{N}(t,t')=r_{N}(2^k \xi^{Rs}_{k,l}t) r_{N}(2^j s^2 \xi^{Rs}_{k,l}t')$ and 
    \begin{align*} D^{m,k}_{s,R}(t,l)f(x)&=\int_{\mathbb{R}^n} \varphi_{m} (2^{(1+\epsilon)k} \xi^{Rs}_{k,l}) r_{N}(2^k \xi^{Rs}_{k,l}t) e^{2\pi\iota x\cdot\xi} \hat{f}(\xi) ~d\xi,~~\text{for}~0\leq m\leq N-1.
\end{align*}
The required $L^p$--estimates of the maximal function associated with \[ \int_{\mathbb{R}} \hat{\psi}(t) e^{2\pi\iota d^{\epsilon}_{k,l}t} D^{m,k}_{s,R}f(x) ~dt\] 
follow by using the arguments as in the case of the maximal function $$\sup_{R>0}\sup_{s\in [s_1,s_2]}|X^k_{N,R,s}f(\cdot)|.$$

Next, consider  
\[E(t,t')=e^{2\pi\iota d^{\epsilon}_{k,l}t} e^{2\pi\iota 2^j(1-a-s^22^{-(1+\epsilon)k}l)t'}\] and  
\begin{align*}
   I^k_{N,R,s}f(x)&=\int_{\mathbb{R}^2} \widehat{\tilde{\psi}}(t)\hat{\psi}(t') E(t,t') T^{Rs}_{N}(t,t')f(x)~dt dt', 
\end{align*}
   \text{where}
   \begin{align*} T^{Rs}_{N}(t,t')f(x)&=\int_{\mathbb{R}^n} s^{R}_{N}(t,t') \varphi\Big(2^{(1+\epsilon)k} \xi^{Rs}_{k,l}\Big) e^{2\pi\iota x\cdot\xi} \hat{f}(\xi) ~d\xi.
   	\end{align*}
We write \[\widehat{K^{R,s}_{t,t'}}(\xi)=2^{j\epsilon N}2^{d(1+\epsilon)N} 2^{k(1+\epsilon)N} s^{R}_{N}(t,t') \varphi(2^{(1+\epsilon)k}\xi^{Rs}_{k,l})\] and note that an integration by parts argument easily gives us that  
\begin{align*}
    |T^{Rs}_{N}(t,t')f(x)|&=2^{-j\epsilon N} 2^{-d(1+\epsilon)N} 2^{-k(1+\epsilon)N} |K^{R,s}_{t,t'}\ast f(x)|\\&\lesssim 2^{-j\epsilon N} 2^{-d(1+\epsilon)N} 2^{-k(1+\epsilon)N} \mathfrak{M}f(x).
\end{align*}
Consequently, we get the required weighted estimates from that of $\mathfrak{M}$.
\qed

\subsection{Proof of Lemma \ref{Tdeltaweightedestimate01}}
	Let $n_0>1$ be such that $\nu+n_0>n/2$. Consider  
	\[T^{\nu}_{t}g(x)=\sum^{n_0}_{k=1}(T^{\nu+k-1}g(x)-T^{\nu+k}g(x))+T^{\nu+n_0}_{t}g(x).\]
	Therefore, 
	\begin{align*}
		\sup_{R>0}\Big(\frac{1}{R\sqrt{2^{1-j}}}\int^{R\sqrt{2^{1-j}}}_{0}|T^{\nu}_{t}g(x)|^{2}~dt\Big)^{1/2}&\leq \sum^{n_0}_{k=1}\sup_{R>0}\Big(\frac{1}{R\sqrt{2^{1-j}}}\int^{R\sqrt{2^{1-j}}}_{0}|\tilde{T}^{\nu+k-1}_{t}g(x)|^{2}~{dt}\Big)^{1/2}\\
		&+T^{\nu+n_0}_{*}g(x),
	\end{align*}
	where $T^{\nu+n_0}_{*}$ denotes the maximal cone multiplier operator with index $\nu+n_0$ and the Fourier multiplier of the operator $\tilde{T}^{\nu+k-1}_{t}$ is given by  $(1-\frac{|\eta'|^2}{t^2\eta^2_n})^{\nu+k-1}_{+}\frac{|\eta'|^2}{t^2\eta^2_n}\varphi(\eta_n)$. 
	
	Since $\nu+n_0>n/2$, we invoke  Proposition $7.2$ of \cite{CHLY}  to deduce that the term corresponding to $T^{\nu+n_0}_{*}$ is pointwise dominated by the Hardy--Littlewood maximal function $M_{HL}$, and hence it is bounded on $L^{2}(\omega_{\alpha,\beta})$ for $\alpha\in [0,n-1), \beta\in [0,1)$.  
	
	Next, for $1\leq k\leq n_0$, it is enough to prove that for $\nu>\max\{\frac{\alpha+\beta}{2}-1,-\frac{1}{2}\}$ the following estimate holds
	\[\Big\Vert \sup_{R>0}\Big(\frac{1}{R\sqrt{2^{1-j}}}\int^{R\sqrt{2^{1-j}}}_{0}|\tilde{T}^{\nu}_{t}g(x)|^{2}~{dt}\Big)^{1/2}\Big\Vert_{L^{2}(\omega_{\alpha,\beta})}\leq C\Vert g\Vert_{L^{2}(\omega_{\alpha,\beta})}.\]
	Let $\Psi\in C^{\infty}_{c}([1/4,1])$ such that \[\sum_{\gamma\geq0}\Psi(2^{\gamma}t)=1.\]
	Therefore, 
	\begin{eqnarray}\label{sum2}
		|\tilde{T}^{\nu}_{t}g(x)|\leq \sum_{\gamma\geq0}2^{-\nu\gamma}|\tilde{T}^{\nu}_{\gamma,t}g(x)|,
	\end{eqnarray}
	where the multiplier of the operator $\tilde{T}^{\nu}_{\gamma,t}, \gamma\geq1$, is given by 
	\begin{eqnarray}\label{multipliersymbol}
		\Psi\Big(2^{\gamma}(1-\frac{|\eta'|^2}{t^{2}\eta^2_n})\Big)2^{\nu\gamma}(1-\frac{|\eta'|^2}{t^2\eta^{2}_n})^{\nu}_{+}\frac{|\eta'|^2}{t^2\eta^{2}_n}\varphi(\eta_n),
	\end{eqnarray}
	and the multiplier of the operator $\tilde{T}^{\nu}_{0,t}$ is \[\Big(1-\sum_{\gamma\geq1}\Psi\Big(2^{\gamma}(1-\frac{|\eta'|^2}{t^{2}\eta^2_n})\Big)\Big)2^{\nu\gamma}(1-\frac{|\eta'|^2}{t^2\eta^{2}_n})^{\nu}_{+}\frac{|\eta'|^2}{t^2\eta^{2}_n}\varphi(\eta_n).\]
	
	It is easy to observe that the maximal function corresponding to the operator $\tilde{T}^{\nu}_{0,t}$ is controlled by $\mathfrak{M}$, that is 
	\[ \tilde{T}^{\nu}_{0,*}g(x)\lesssim \mathfrak{M}g(x).\]
	Therefore, the desired $L^{2}(\omega_{\alpha,\beta})$ estimates for this term follows whenever  $\alpha\in [0,n-1), \beta\in [0,1)$. 
	
	On the other hand notice that 
	\[ \sup_{R>0}\Big(\frac{1}{R\sqrt{2^{1-j}}}\int^{R\sqrt{2^{1-j}}}_{0}|\tilde{T}^{\nu}_{\gamma,t}g(x)|^{2}~{dt}\Big)^{1/2}\leq \Big(\int^{\infty}_{0} |\tilde{T}^{\nu}_{\gamma,t}g(x)|^{2}\frac{dt}{t}\Big)^{1/2}.\]
	In this case, invoking Lemma \ref{Tpsiweightedestimate1} yields that 
	\begin{eqnarray}\label{weightedinequality}
		\Vert \Big(\int^{\infty}_{0} |\tilde{T}^{\nu}_{\gamma,t}g(x)|^{2}\frac{dt}{t}\Big)^{1/2}\Vert_{L^{2}(\omega_{\alpha,\beta})}\leq C A_{2^{-\gamma}}(\alpha,\beta)\Vert g\Vert_{L^{2}(\omega_{\alpha,\beta})}.  
	\end{eqnarray}
	Hence, by using the estimate \eqref{weightedinequality} in \eqref{sum2}, we complete the proof of Lemma~\ref{Tdeltaweightedestimate01}.  
	\qed
	
Note that with this we have completed the proof of desired weighted $L^2\times L^2 \to L^1$-estimates of the operators $\mathcal{T}^{\lambda}_{*,j}$ for $j \geq 2$. Therefore, in order to complete the proof of Theorem~\ref{mainresult}, we are left with proving the required weighted $L^2\times L^2 \to L^1$-estimates for the operator $\mathcal{T}^{\lambda}_{*,1}$. We take up this work  in the following subsection.  

\subsection{Boundedness of $\mathcal{T}^{\lambda}_{*,1}(f,g)$}\label{T0}
Recall that the multiplier of the bilinear operator $\mathcal{T}^{\lambda}_{R,1}$ is given by 
\[m^{\lambda,1}(\frac{\xi'}{R\xi_n},\frac{\eta'}{R\eta_n})=\psi_{1}(\frac{|\xi'|^2}{R^2\xi^2_{n}}) \Big(1-\frac{|\xi'|^2}{R^2\xi^2_n}-\frac{|\eta'|^2}{R^2\eta^2_n}\Big)^{\lambda}_{+}\Phi(\xi_n,\eta_n),\]
where $\Phi(\xi_n,\eta_n)=\varphi(\xi_n)\varphi(\eta_n)$.

We begin with a similar decomposition of the multiplier $m^{\lambda,1}(\frac{\xi'}{R\xi_n},\frac{\eta'}{R\eta_n})$ as we did in the previous section. However, this time we need to perform the decomposition with respect to the variable $\eta$. Consider 
\begin{eqnarray*}
	&&m^{\lambda,1}(\frac{\xi'}{R\xi_n},\frac{\eta'}{R\eta_n})\\
	&&=\psi_{1}(\frac{|\xi'|^2}{R^2\xi^2_{n}})\sum_{j\geq2}\psi\Big(2^{j}(1-\frac{|\eta'|^2}{R^2\eta^2_n})\Big)(1-\frac{|\eta'|^2}{R^2\eta^2_n})^{\lambda}_{+} \Big(1-\frac{|\xi'|^2}{R^2\xi^{2}_n}(1-\frac{|\eta'|^2}{R^2\eta^2_n})^{-1}\Big)^{\lambda}_{+}\Phi(\xi_n,\eta_n)\\
	&&+\psi_{1}(\frac{|\xi'|^2}{R^2\xi^2_{n}})\psi_{1}(\frac{|\eta'|^2}{R^2\eta^2_{n}}) (1-\frac{|\xi'|^2}{R^2\xi^2_{n}}-\frac{|\eta'|^2}{R^2\eta^2_{n}})^{\lambda}_{+}\Phi(\xi_n,\eta_n).
\end{eqnarray*}
Further, we split $\psi_1=\psi^1_1+\psi^2_1$ such that $\text{supp}(\psi^1_1)\subset [0,\frac{3}{16}]$ and $\text{supp}(\psi^2_1)\subset [\frac{3}{32},\frac{3}{4}]$. Then, we have 

\begin{align*}
	m^{\lambda,1}(\frac{\xi'}{R\xi_n},\frac{\eta'}{R\eta_n})
	&=\sum_{j\geq 2} m^{\lambda,1}_{j}(\frac{\xi'}{R\xi_n},\frac{\eta'}{R\eta_n})+m^{\lambda,1}_{1,1}(\frac{\xi'}{R\xi_n},\frac{\eta'}{R\eta_n})+m^{\lambda,1}_{1,2}(\frac{\xi'}{R\xi_n},\frac{\eta'}{R\eta_n}),
\end{align*}
where
\begin{align*}
	m^{\lambda,1}_{j}(\frac{\xi'}{R\xi_n},\frac{\eta'}{R\eta_n})
	&=\psi_{1}(\frac{|\xi'|^2}{R^2\xi^2_{n}})\Phi(\xi_n,\eta_n)\psi\Big(2^{j}(1-\frac{|\eta'|^2}{R^2\eta^2_n})\Big)(1-\frac{|\eta'|^2}{R^2\eta^2_n})^{\lambda}_{+} \Big(1-\frac{|\xi'|^2}{R^2\xi^{2}_n}(1-\frac{|\eta'|^2}{R^2\eta^2_n})^{-1}\Big)^{\lambda}_{+},~j\geq 2
\end{align*}
and
\begin{align*} 
	m^{\lambda,1}_{1,i}(\frac{\xi'}{R\xi_n},\frac{\eta'}{R\eta_n})&=\Phi(\xi_n,\eta_n)\psi_{1}(\frac{|\xi'|^2}{R^2\xi^2_{n}})\psi^i_{1}(\frac{|\eta'|^2}{R^2\eta^2_{n}}) (1-\frac{|\xi'|^2}{R^2\xi^2_{n}}-\frac{|\eta'|^2}{R^2\eta^2_{n}})^{\lambda}_{+}, ~i=1,2.	
\end{align*}
Observe that on the support of $m^{\lambda,1}_{1,1}$, the quantity $\frac{|\xi'|^2}{R^2 \xi^2_n}+\frac{|\eta'|^2}{R^2 \eta^2_n}\leq \frac{3}{4}+\frac{3}{16}<1,$ which implies that   $m^{\lambda,1}_{1,1}$ is a smooth function. Therefore, 
\[\sup_{R>0}\Big|\int_{\mathbb{R}^{2n}}m^{\lambda,1}_{1,1}(\frac{\xi'}{R\xi_n},\frac{\eta'}{R\eta_n}) \hat{f}(\xi)\hat{g}(\eta) e^{2\pi\iota x\cdot(\xi+\eta)} d\xi d\eta   \Big|\lesssim \mathfrak{M}f(x) \mathfrak{M}g(x).\]

Invoking the weighted boundedness of the operator $\mathfrak{M}$ we get that the operator associated with $m^{\lambda,1}_{1,1}$ is bounded from $L^2(\omega_{\alpha_1,\beta_1})\times L^{2}(\omega_{\alpha_2,\beta_2})\to L^1(\omega)$ for $\alpha_i\in [0,n-1)$, $\beta_i\in [0,1), i=1,2$.
\subsection*{Estimates of the operator corresponding to $\sum_{j\geq2} m^{\lambda,1}_{j}$}
Let $K>2$ be such that $\frac{|\xi'|}{\xi_n}\leq \frac{R}{8}$ for all $(\xi',\xi_n)\in \text{supp}(m^{\lambda,1}_{j})$ for $j\geq K$.  Note that in this case $\psi_{1}(\frac{|\xi'|^2}{R^2\xi^2_n})\equiv 1$. 
Therefore, we get that  
\begin{align*}
	\mathcal{T}^{\lambda,1}_{K,R}(f,g)(x)
	&=\sum_{j\geq K} \int_{\mathbb{R}^{2n}}m^{\lambda,1}_{j}(\frac{\xi'}{R\xi_n},\frac{\eta'}{R\eta_n}) \hat{f}(\xi)\hat{g}(\eta) e^{2\pi\iota x\cdot (\xi+\eta)} d\xi d\eta.
\end{align*}
Hence, the estimate of $\sup_{R>0}|\mathcal{T}^{\lambda,1}_{K,R}(f,g)|$ follows from the boundedness of the operator $\mathcal{T}^{\lambda}_{*,j}, j\geq2$ with the roles of $f$ and $g$ interchanged.

Now it remains to prove boundedness of the operator corresponding to the multiplier $m^{\lambda,1}_{j}$ for $2\leq j<K.$
This case is dealt with in a similar manner as it was done for the multipliers $m^{\lambda,j}, j\geq 2,$ however with the roles of $\xi$ and $\eta$ interchanged. We will briefly discuss it to point out the differences. However, we omit the details to avoid repetition. Note that, as compared to the previous case, here we get an extra term $\psi_{1}(\frac{|\xi'|^2}{R^2\xi^2_n})$ with the multiplier. We can club it with $\hat{f}(\xi)$. Further, the Cauchy-Schwarz inequality reduces the matter to proving suitable weighted $L^2$-estimates  of the two sublinear operators acting on $f$ and $g$ separately. The boundedness of the operator corresponding to $g$ follows from Proposition \ref{prop5.4}. 
In order to deal with the operator corresponding to $f$, we consider  
\[B^{\psi_{1}}_{1,R}f(x):=\int_{\mathbb{R}^n}\varphi(\xi_n)\psi_{1}(\frac{|\xi'|^2}{R^2\xi^2_n}) \hat{f}(\xi) e^{2\pi\iota x\cdot \xi}~d\xi.\]
The following result holds for the operator as above.
\begin{lemma}
	Let $n\geq3$ and $\alpha_1\in[0,n-1)$ and $\beta_1\in [0,1)$. Then for $\nu>\max\{\frac{\alpha_1+\beta_1}{2}-1,-\frac{1}{2}\}$, we have 
	\begin{align*}\|\sup_{R>0}\Big(\frac{1}{R_j}\int^{R_j}_{0}|B^{\psi_{1}}_{1,R}T^{\nu}_{t}f(x)|^2~dt\Big)^{1/2}\|_{L^{2}(\omega_{\alpha_1,\beta_1})}\lesssim \|f\|_{L^{2}(\omega_{\alpha_1,\beta_1})}.
		\end{align*}
\end{lemma}
\begin{proof} Let $n_0$ be a large enough number such that $\nu+n_0>n/2$. Consider 
	\begin{align*}
		\sup_{R>0}\Big(\frac{1}{R_j}\int^{R_j}_{0}|B^{\psi_{1}}_{1,R}T^{\nu}_{t}f(x)|^2~dt\Big)^{1/2}&\leq \sum^{n_0}_{k=1} \sup_{R>0}\Big(\frac{1}{R_j}\int^{R_j}_{0}|B^{\psi_{1}}_{1,R}(T^{\nu+k}_{t}-T^{\nu+k-1}_{t})f(x)|^2~dt\Big)^{1/2}\\
		&+ \sup_{R>0}\Big(\frac{1}{R_j}\int^{R_j}_{0}|B^{\psi_{1}}_{1,R}T^{\nu+n_0}_{t}f(x)|^2~dt\Big)^{1/2}.
	\end{align*}
	Since $\sup_{t>0}|B^{\psi_{1}}_{1,R}T^{\nu+n_0}_{t}f(x)|\lesssim \mathfrak{M}(T^{\nu+n_0}_{*}f)(x)\lesssim \mathfrak{M}\circ M_{HL}f(x)$ for a.e. $x$, we get the desired  $L^2(\omega_{\alpha_1,\beta_1})$ estimate of this term for $\alpha_1\in [0,n-1)$ and $\beta_1\in [0,1)$.
	On the other hand, 
	\begin{align*}
		|(T^{\nu+1}_t-T^{\nu}_t)f(x)|\lesssim \sum_{\gamma\geq0} 2^{-\nu\gamma}|\tilde{T}^{\nu}_{\gamma,t}f(x)|,
	\end{align*}
	where the multiplier corresponding to $\tilde{T}^{\nu}_{\gamma,t}$ is given by \eqref{multipliersymbol}. Next, observe that $\tilde{T}^{\nu}_{0,*}f(x)\lesssim \mathfrak{M}f(x),$ therefore, 
	\begin{align*}\sup_{R>0}|B^{\psi_0}_{1,R}\tilde{T}^{\nu}_{0,*}f(x)|&\lesssim \mathfrak{M}\circ \mathfrak{M}f(x).\end{align*}
	Consequently, the desired $L^{2}(\omega_{\alpha_1,\beta_1})$ estimate follows for this term as well. 
	Next, we have 
	\begin{align*}
		\sum_{\gamma\geq1} 2^{-\nu\gamma} \sup_{R>0}\Big(\frac{1}{R_j}\int^{R_j}_0|B^{\psi_0}_{1,R}(\tilde{T}^{\nu}_{\gamma,t}f(x))|^2~dt\Big)^{1/2}\leq \sum_{\gamma\geq1} 2^{-\nu\gamma} \Big(\int^{\infty}_{0}|\mathfrak{M}(\tilde{T}^{\nu}_{\gamma,t}f)(x)|^2 \frac{dt}{t}\Big)^{1/2}.
	\end{align*}
	Therefore, 
	\begin{align*}
		\Vert \Big(\int^{\infty}_{0}|\mathfrak{M}(\tilde{T}^{\nu}_{\gamma,t}f)(x)|^2 \frac{dt}{t}\Big)^{1/2}\Vert_{L^{2}(\omega_{\alpha_1,\beta_1})}&= \Big(\int^{\infty}_{0}\int_{\mathbb{R}^n}|\mathfrak{M}(\tilde{T}^{\nu}_{\gamma,t}f)(x)|^2 \omega_{\alpha_1,\beta_1}(x) ~dx \frac{dt}{t}\Big)^{1/2}\\
		&\lesssim \Big(\int_{\mathbb{R}^n}\int^{\infty}_{0} |\tilde{T}^{\nu}_{\gamma,t}f(x)|^2 \frac{dt}{t} \omega_{\alpha_1,\beta_1}(x)~dx\Big)^{1/2}\\
		&\lesssim A_{2^{-\gamma}}(\alpha_1,\beta_1) \Vert f\Vert_{L^{2}(\omega_{\alpha_1,\beta_1})}.
	\end{align*}
	The final inequality in the equation above follows from equation \eqref{weightedinequality}. 
	
	It remains to prove boundedness of the operator with multiplier $m^{\lambda,1}_{1,2}$.
	As earlier using identity~\eqref{SteinWeiss} we can split the multiplier as follows
	\begin{eqnarray*}
		m^{\lambda,1}_{1,2}(\frac{\xi'}{R\xi_n},\frac{\eta'}{R\eta_n})=\Phi(\xi_n,\eta_n)\psi_{1}(\frac{|\xi'|^2}{R^2\xi^2_n})\psi^{2}_{0}(\frac{|\eta'|^2}{R^2\eta^2_n}) R^{-2\alpha}\int^{uR}_{0}(R^2\phi_{R}(\eta)-t^2)^{\mu-1}_{+} t^{2\nu+1} (1-\frac{|\xi'|^2}{t^2\xi^2_n})^{\nu}_{+}~dt, 
	\end{eqnarray*}
	where  $u=\sqrt{\frac{29}{32}}$, $\lambda=\mu+\nu$ with $\mu>0$ and $\nu>-1$.
	The Cauchy--Schwarz inequality and a change of variable $t\to Rt$ reduce the problem to proving suitable weighted $L^2$-estimates of the operator 
	\begin{align*}\sup_{R>0}\int^u_{0}|H_{\mu}^{R,t}g(\cdot)|^2 ~dt,
		\end{align*}
	where \[H_{\mu}^{R,t}g(x)=\int_{\mathbb{R}^n} \varphi(\eta_n)\psi^{2}_{0}(\frac{|\eta'|^2}{R^2\eta^2_n}) (1-t^2-\frac{|\eta'|^2}{R^2\eta^2_n})^{\mu-1}_{+} \hat{g}(\eta) e^{2\pi\iota x\cdot \eta}~d\eta.\]
	The operator $H_{\mu}^{R,t}$ is similar to $S^{R,t}_{j,\mu}$. Therefore, similar arguments prove that  \[\sup_{R>0}\int^u_{0}|H_{\mu}^{R,t}g(\cdot)|^2 ~dt\] is bounded on $L^{2}(\omega_{\alpha_2,\beta_2})$ for $\alpha_2\in [0,n-1), \beta_2\in [0,1)$ and $\mu>\max\{\frac{\alpha_2+\beta_2}{2},\frac{1}{2}\}$. We omit the details to avoid repetition. This completes the proof. 
\end{proof}

\section{Appendix}\label{Appendix}
We give a direct proof of $L^2(\R^n)\times L^2(\R^n)\to L^1(\R^n)$--boundedness of the maximal operator $\mathcal{T}^{\lambda}_{*}$ for $\lambda>0$. 

We begin with the same decomposition of the operator $\mathcal{T}^{\lambda}_{*}$ as done in  equation~\eqref{operatorT_j}. For $j\geq2$, we have that 
        \begin{align*}|\mathcal{T}^{\lambda}_{R,j}(f,g)(x)|\leq c_{\lambda} 2^{-j/4} \Big( \int^{\sqrt{2^{1-j}}}_{0} |S^{R,t}_{j,\mu}f(x) t^{2\nu+1}|^2 ~dt\Big)^{1/2} \Big(\frac{1}{R\sqrt{2^{1-j}}}\int^{R\sqrt{2^{1-j}}}_{0}|T^{\nu}_{t}f(x)|^2 dt\Big)^{1/2}.
        	\end{align*}
		
Therefore, the desired $L^2(\R^n)\times L^2(\R^n)\to L^1(\R^n)$--boundedness of the maximal operator $\mathcal{T}^{\lambda}_{*}$ for $\lambda>0$ would follow from $L^2$--estimates of the maximal operators corresponding to the terms on the right hand side of the inequality above.  Moreover, the term corresponding to $j=1$ in the decomposition of $ \mathcal{T}^{\lambda}_{R,j}(f,g)(x)$ is dealt with in a similar fashion as in Subsection \ref{T0}. Since the proof is repetitive with arguments similar to the case of $j\geq 2$, we skip the details. 

We give the proof of estimates for terms with $j\geq 2$. In this regard, we have the  following results. 
		\begin{lemma}\label{fullsquarefunction}
			For  $n\geq3$ and $\nu>-1/2$, we have 
			\[ \Big\Vert \sup_{R>0}\Big(\frac{1}{R\sqrt{2^{1-j}}}\int^{R\sqrt{2^{1-j}}}_{0}|T^{\nu}_{t}f(x)|^2 dt\Big)^{1/2}\Big\Vert_{L^2(\R^n)}\lesssim \Vert f\Vert_{L^{2}(\R^n)}.\]
		\end{lemma}
	
		\begin{lemma}\label{lemma2appendix}
			For $n\geq3$
              and $\mu>1/2$, we have 
            \[\Big\Vert \sup_{R>0}\Big(\int^{\sqrt{2^{1-j}}}_{0}|S^{R,t}_{j,\mu}f(x)t^{2\nu+1}|^2~dt\Big)^{1/2} \Big\Vert_{L^{2}(\R^n)}\lesssim   2^{j/4} 2^{-j(\mu+\nu)}  \Vert f\Vert_{L^{2}(\R^n)}.\]
		\end{lemma}
			\subsection*{Proof of Lemma \ref{fullsquarefunction}}
			Note that
			\begin{align*}
				\sup_{R>0}\Big(\frac{1}{R\sqrt{2^{1-j}}}\int^{R\sqrt{2^{1-j}}}_{0}|T^{\nu}_{t}f(x)|^2 dt\Big)^{1/2}&=\sup_{R>0}\Big( \frac{1}{R} \int^R_{0} |T^{\nu}_{t}f(x)|^2 ~dt\Big)^{1/2}\\
				&:=M^{\nu}f(x).
			\end{align*}
			The $L^2$--boundedness of $M^{\nu}f$ is deduced by the corresponding estimates of the square function associated with cone multipliers. Let $n\geq3$ and $\nu>-1/2$. Consider the square function associated with the linear cone multipliers  
			\[\mathcal{G}^{\nu}f(x)=\Big( \int^{\infty}_0 \Big| T^{\nu+1}_t f(x)-T^{\nu}_{t}f(x)\Big|^2 \frac{dt}{t}\Big)^{1/2},\]
			where $T^{\nu}_{t}f(x)=\int_{\R^n} \varphi(\eta_n)(1-\frac{|\eta'|^2}{t^2 \eta^2_n})^{\nu}_{+} \hat{f}(\eta) e^{2\pi\iota x\cdot\eta}~d\eta$. 
			
			The following $L^2$--estimate holds for the square function. 
			\begin{proposition}\label{squarefunction}
					Let $n\geq3$ and $\nu>-1/2$. Then 
					\[\Vert\mathcal{G}^{\nu}f\Vert_{L^2(\R^n)}\lesssim \Vert f\Vert_{L^{2}(\R^n)}. \]
				\end{proposition}
The proof of this proposition follows from Plancherel's theorem. Indeed, consider  
					\begin{align*}
						\Vert \mathcal{G}^{\nu}f\Vert^2_{L^2(\R^n)}&=\int_{\R^n} \int^{\infty}_0 \Big| T^{\nu+1}_t f(x)-T^{\nu}_{t}f(x)\Big|^2 \frac{dt}{t} ~dx\\
						\text{(by Plancherel's theorem)}&=\int^{\infty}_0 \int_{\R^n} |\hat{f}(\eta) \varphi(\eta_n)|^2 \left(1-\frac{|\eta'|^2}{t^2\eta^2_n}\right)^{2\nu}_{+} \frac{|\eta'|^4}{(t\eta_n)^4} d\eta ~\frac{dt}{t}\\
						\text{(by Fubini's theorem)}&= \int_{\R^n} |\hat{f}(\eta) \varphi(\eta_n)|^2 \frac{|\eta'|^4}{\eta^4_n} \int^{\infty}_{|\eta'|/\eta_n}  \left(1-\frac{|\eta'|^2}{t^2\eta^2_n}\right)^{2\nu}~\frac{dt}{t^5}~d\eta\\
						\text{(by change of variable)}&= \int_{\R^n} |\hat{f}(\eta) \varphi(\eta_n)|^2 \int^{\infty}_1 t^{-4\nu-5} (t^2-1)^{2\nu} ~dt ~d\eta\\
						&=\int_{\R^n} |\hat{f}(\eta) \varphi(\eta_n)|^2 \int^{\infty}_1 t^{-2\nu-3} (t-1)^{2\nu}~dt~d\eta\\
						(\text{since} ~\nu>-1/2) &\lesssim \int_{\R^n} |\hat{f}(\eta) \varphi(\eta_n)|^2 ~d\eta\\
						&\lesssim \Vert f\Vert^2_{L^2(\R^n)}.
					\end{align*}
					This completes the proof of Proposition~\ref{squarefunction}.

				Next, observe that 
				\begin{align*}M^{\nu}f(x)\leq \mathcal{G}^{\nu}f(x)+M^{\nu+1}f(x).
					\end{align*}
					Further, employing the above process for $k$-times where $k+\nu>n/2$, we get that
					\begin{eqnarray}\label{relation}
						M^{\nu}f(x)\leq \mathcal{G}^{\nu}f(x)+\mathcal{G}^{\nu+1}f(x)+\mathcal{G}^{\nu+2}f(x)+\cdots+\mathcal{G}^{\nu+k-1}f(x)+ M^{\nu+k}f(x).
					\end{eqnarray}
					Observe that $ M^{\nu+k}f(x)\leq T^{\nu+k}_{*}f(x)\lesssim M_{HL}f(x)$ for $\nu+k>n/2$ (see \cite{CHLY}). Therefore, the estimate in Lemma~\ref{fullsquarefunction} follows by invoking $L^2$--estimate of the square function $\mathcal{G}^{\nu}f$ from Proposition~\ref{squarefunction}. This completes the proof of Lemma~\ref{fullsquarefunction}. \qed
					
        \subsection*{Proof of Lemma~\ref{lemma2appendix}}
            
	We give a sketch of the proof of Lemma \ref{lemma2appendix}.  We require the following two results along with the decompositions of the operator from Section \ref{mainproof}.
        \begin{lemma}\label{lem1}
            Let $n\geq3$ and $0<\delta\ll 1$. Then 
            \[\Big\Vert \Big(\int^{\infty}_0 |T^{\Psi}_{\delta,t}f(\cdot)|^2 \frac{dt}{t}\Big)^{1/2}\Big\Vert_{L^2(\R^n)} \lesssim \delta^{1/2}\Vert f\Vert_{L^2(\R^n)}.\]
        \end{lemma}
        The proof of Lemma \ref{lem1} follows by using Plancherel's theorem. Further, Lemma \ref{lem1} and the Sobolev embedding theorem yield the following $L^2$--estimate of the corresponding maximal function. 
       \begin{lemma}\label{lem2}
            Let $n\geq3$ and $0<\delta\ll 1$. Then 
            \[\Big\Vert \sup_{t>0}|T^{\Psi}_{\delta,t}f(\cdot)|\Big\Vert_{L^2(\R^n)} \lesssim \Vert f\Vert_{L^2(\R^n)}.\]
        \end{lemma}
 
	We follow the line of arguments as were used in proving Lemma~\ref{SRtjbeta} in Subsection~\ref{case2}. As in there, we consider the cases of $t<\sqrt{2^{-1-j}}$ and $t\geq \sqrt{2^{-1-j}}$ separately. 
	
	Recall that in {\bf Case I}, the operator $S^{R,t}_{j,\mu}$ takes the form \eqref{Appen1}
	\[S^{R,t}_{j,\gamma+1}f(x)=2^{j\rho}\sum_{k\geq0} \frac{\Gamma(\rho+k)}{\Gamma(\rho)k!} (2^{j}t^{2})^{k} T_{2^{-j},R}^{\psi^{k}}f(x),~~-1<\gamma<0.\]
	We invoke Lemma \ref{lem2} to get that  
	\begin{align*}\Vert \sup_{R>0}|S^{R,t}_{j,\gamma+1}f|\Vert_{L^{2}(\R^n)}&\lesssim 2^{j\rho} \sum_{k\geq0} \frac{\Gamma(\rho+k)}{\Gamma(\rho)k!} 2^{-(1+\epsilon_{0})k} k^{N} 2^{N+k} \Vert f\Vert_{L^{2}(\R^n)}.
		\end{align*}
	Consequently, we get that 
	\begin{align*}
		\Big\Vert \sup_{R>0}\Big(\int^{\sqrt{2^{1-j}}}_{0}|S^{R,t}_{j,\mu}f(x)t^{2\nu+1}|^2~dt\Big)^{1/2} \Big\Vert_{L^{2}(\R^n)}&\lesssim 2^{j\rho} 2^{-j(4\nu+3)/4} \Vert f\Vert_{L^{2}(\R^n)}\\
        &= 2^{-j(\mu+\nu)} 2^{j/4} \Vert f\Vert_{L^{2}(\R^n)}.
	\end{align*}
    The remaining cases $\gamma\geq0$ can be completed using the same arguments as in {\bf Case I} of Subsection \ref{case2} along with Lemma \ref{lem2}.
    
    Next, for the {\bf Case II}, that is when $t^2\in [2^{-j-1-\epsilon_{0}}, 2^{-j+1}]$, in view of the equation~\eqref{Appen2}, we need to prove the following $L^2(\R^d)$--estimate. Note that this is an analogue of Lemma~\ref{prop5.3}.  

	\begin{proposition}
		Let  $2\leq j<M'$. Then for $\mu>1/2$, we have 
		\[\Big\Vert \sup_{R>0}\Big(\int^{s_2}_{s_1}|S'^{R,s}_{j,\mu}f(x)|^2~ds\Big)^{1/2}\Big\Vert_{L^{2}(\R^n)}\leq C_{M',n} \Vert f\Vert_{L^{2}(\R^n)}.\]
	\end{proposition}
	
	\begin{proof}
		Note that we have 
		\[S'^{R,s}_{j,\mu}f(x)=T^{\psi}_{2^{-j},R}\circ T^{\mu-1}_{sR}f(x),\]
		with the following estimate (by equation~\eqref{Maximalfunction1})
		
        \begin{eqnarray*}
        \sup_{R>0}|T^{\psi}_{2^{-j},R}f(x)|\lesssim C_j \mathfrak{M}f(x),~\text{a.e.}~x.
        \end{eqnarray*}
As earlier, recall that for a large integer $n_0$ satisfying $\mu-1+n_0>n/2$, we have 
		\begin{eqnarray*}
			&&\sup_{R>0}\Big(\int^{s_2}_{s_1}|S'^{R,s}_{j,\mu}f(x)|^2~ds\Big)^{1/2}\\
			&&\lesssim 2^{M'n}\Big(\int^{\infty}_{0}|\mathfrak{M}\circ \tilde{T}^{\mu-1}_{s}f(x)|^2~\frac{ds}{s}\Big)^{1/2}+ 2^{M'n} \mathfrak{M}\circ T^{\mu-1+n_0}_{*}f(x).
		\end{eqnarray*}
        
		Therefore, using $L^2$--boundedness of the operators $\mathfrak{M}$, $M_{HL}$,
		and Lemma \ref{lem1}, we get the desired  $L^{2}$--estimate of the operator under consideration  for $\mu>1/2$.
	\end{proof}
	
	Finally, in order to complete the proof, we need to prove the following analogue of Lemma~\ref{prop5.4} for large values of the parameter $j$. 
	\begin{proposition}
		Let $n\geq3$, $\mu>\frac{1}{2}$ and $j\geq M'$ for some large $M'\in\mathbb{N}$. Then we have 
		\[\Big\Vert \sup_{R>0}\Big(\int^{s_2}_{s_1}|S'^{R,s}_{j,\mu}f(x)|^2~ds\Big)^{1/2} \Big\Vert_{L^{2}(\R^n)}\lesssim 2^{-j(\mu-1)}  2^{-j/2}  \Vert f\Vert_{L^{2}(\R^n)}.\]
	\end{proposition}

We observe that by following the proof of Lemma~\ref{prop5.4}, we only need to prove the desired $L^2$-estimates of the operators $U^{\phi_q,j}_{Rs, R}(k,l)$ and $P^{j,k}_{Rs,R}(t,l)$. Further, the $L^2$--boundedness of the operator $U^{\phi_q,j}_{Rs, R}(k,l)f$ follows by using the inequality \eqref{Vphi} and Lemma \ref{lem1}. Moreover, the $L^2$--boundedness of the operator  
$P^{j,k}_{Rs,R}(t,l)f$ follows from $L^2$--boundedness of the maximal function $\mathfrak{M}f$. This is similar as in Subsection \ref{Pjk}. This completes the proof.
\qed

\section*{Acknowledgment}
Saurabh Shrivastava acknowledges the financial support from Science and Engineering Research Board, Govt. of India, under the scheme of Core Research Grant, file no. CRG/2021/000230. Kalachand Shuin acknowledges the support of Inspire Faculty Award of Department of Science and Technology, Govt. of India, Registration no. IFA23-MA191.
\bibliography{biblio}
\end{document}